\documentclass[12pt,reqno]{amsart}
\usepackage[utf8]{inputenc}
\usepackage[colorinlistoftodos]{todonotes}
\usepackage{epsfig}
\usepackage{graphicx}
\usepackage{amsmath}
\usepackage{caption}
\usepackage{breqn}
\usepackage{bbm}
\usepackage{filecontents}
\usepackage{lipsum}

\usepackage{enumerate}
\usepackage{enumitem}
\usepackage{graphicx}
\usepackage{float}
\usepackage{bbm}
\usepackage{amssymb}
\usepackage{comment}
\usepackage{hyperref}
\usepackage{cleveref}
\usepackage{subcaption}
\date{\vspace{-5ex}}
\usepackage{amsthm, amsfonts, mathrsfs}
\graphicspath{./figures}

\usepackage{listings}
\usepackage{color}
\usepackage[lined,boxed,commentsnumbered]{algorithm2e}
\usepackage{esint}

\def\R{\mathbb{R}}

\def\bP{{\mathbf{P}}}

\def\newmap{{single-slice matching~}}
\def\Mnewmap{{matrix-slice matching~}}
\def\Pmatching{{{matrix-slice matching}}}
\def\opname{{slice-matching operator}}
\def\W{{\mathcal{W}}}

\usepackage[%
  sortcites,
  backend=biber,
  defernumbers=true,
  sorting=nyt,           % sort by name, year, title
  style = numeric,
  natbib=true,
  firstinits=true,
  maxbibnames=99,
]{biblatex}
\newtheorem{theorem}{Theorem}[section]
\newtheorem{proposition}[theorem]{Proposition}
\newtheorem{lemma}[theorem]{Lemma}
\newtheorem{corollary}[theorem]{Corollary}
\newtheorem{definition}[theorem]{Definition}
\newtheorem{remark}[theorem]{Remark}

\newtheorem{example}{Example}[section]

\newtheoremstyle{discussionstyle}
  {10pt} % Space above
  {10pt} % Space below
  {}     % Body font
  {}     % Indent amount
  {\bfseries} % Theorem head font
  {.}    % Punctuation after theorem head
  { }    % Space after theorem head
  {\thmname{#1}\thmnumber{ #2}\thmnote{ (#3)}}

\theoremstyle{discussionstyle}
\newtheorem*{discussion}{Discussion}

\newcommand{\grad}{\triangledown}

\newcommand{\bx}{{{x}}}

\newcommand{\bth}{{{\theta}}}

\newcommand{\bb}{b}
\newcommand{\id}{\operatorname{id}}

\newcommand{\smop}{{\mathcal{U}}}
\newcommand{\affset}{{\mathcal{S}}}

\DeclareMathOperator*{\argmin}{arg\,min}

\title{Approximation properties of slice-matching operators}

\author{Shiying Li}\thanks{Department of Mathematics, University of North Carolina at Chapel Hill, NC, USA. \{shiyl,cmoosm\}$@$unc.edu}
\author{Caroline Moosm\"{u}ller}

\date{}
\begin{document}

%\keywords{Optimal Transport, ...}
%\subjclass[2020]{49Q22, 60D05, 68T10 (needs update)}
%%%%%%%%%%%%%%%%%%%%%%%%%%%%%%%%%%%%%%%%%%%
%%%%%%%%%%%%%%%%%%%%%%%%%%%%%%%%%%%%%%%%%%%
%%%%%%    Abstract                 %%%%%%%%
%%%%%%%%%%%%%%%%%%%%%%%%%%%%%%%%%%%%%%%%%%%

\begin{abstract}
Iterative slice-matching procedures are efficient schemes for transferring a source measure to a target measure, especially in high dimensions. These schemes have been successfully used in applications such as color transfer and shape retrieval, and are guaranteed to converge under regularity assumptions.  In this paper, we explore approximation properties related to a single step of such iterative schemes by examining an associated slice-matching operator, depending on a source measure, a target measure, and slicing directions. In particular, we demonstrate an invariance property with respect to the source measure, an equivariance property with respect to the target measure, and Lipschitz continuity concerning the slicing directions. We furthermore establish error bounds corresponding to approximating the target measure by one step of the slice-matching scheme and characterize situations in which the slice-matching operator recovers the optimal transport map between two measures. We also investigate connections to affine registration problems with respect to (sliced) Wasserstein distances. These connections can be also be viewed as extensions to the invariance and equivariance properties of the slice-matching operator and illustrate the extent to which slice-matching schemes incorporate affine effects.   
\end{abstract}

\maketitle

\setcounter{tocdepth}{1}
\tableofcontents

\section{Introduction}

Optimal transport and Wasserstein distances play a crucial role in machine learning and related applications \cite{arjovsky2017wasserstein,rubner2000earth,kolouri2017optimal,peyre19}. 
For example, these methods have gained prominence in generative modeling \cite{bousquet2017optimal}, aiming to find transport maps or their approximations \cite{baptista2023ApproxFramework, lambert2022variational} that align posterior and prior distributions. A basic problem can be described as follows: Given a random variable $X\sim \sigma$ where $\sigma$ is a prior probability distribution and a target probability distribution $\mu$ of interest,  find  a transformation $T$ such that $T(X)\sim \mu$.

Transforming one probability distribution into another is a fundamental problem that also has implications for flow matching \cite{lipman2022flow}, probability flows \cite{song2021scorebased,chen2019neural}, particle evolutions \cite{taghvaei2016optimal}, and the general learning of underlying distributions in complex data sets 
\cite{arjovsky2017wasserstein, kobyzev2020normalizing}.

 While  Wasserstein distances have proven successful in modeling probability distributions \cite{arjovsky2017wasserstein,lambert2022variational}, their computational expense, especially in high-dimensional scenarios, necessitates more efficient approaches. The computation can be intensive for large-scale problems; specifically, the cost of calculation  via linear programs comes with a complexity of $O(m^3\log(m))$, while the Sinkhorn version \cite{cuturi-2013} provides a faster approximation at $O(m^2\log(m))$, where $m$ is the number of particles  used to approximate a measure.  Sliced-Wasserstein-based generative models \cite{kolouri2018slicedwasserstein, bonet2022efficient} provide scalable alternatives.  Our study focuses on one such model, namely, slice-matching schemes \cite{pitie2007automated,bonneel2015sliced}.  These schemes utilize projections onto lines (slices) as well as the computational benefit of one-dimensional optimal transport to offer effective approximations of target distributions. Furthermore, they are closely aligned with the broader context of normalizing flows \cite{papamakarios2019neural, kobyzev2020normalizing} and variational inference \cite{lambert2022variational}. Beyond their computational advantages, these schemes also demonstrate promising convergence results, as shown in \cite{bonnotte13thesis,pitie2007automated,li2023measure}.

While convergence analysis of iterative slice-matching schemes is the primary focus of \cite{li2023measure}, the present manuscript  aims at establishing a comprehensive understanding of one step of such schemes, focusing on  recovery and stability properties.  We explore structural relationships between measures and address two closely related questions: (A) when can closed-form formulations serve as robust alternatives to optimal transport maps? (B) how effectively do slice-matching approximations represent the target measure? Additionally, our study examines the ability to handle transformations like shifts and scalings in the initial step, demonstrated through the analysis of some basic registration problems. Registration problems are not the main focus of this paper; we rather use them to understand affine effects of slice-matching maps. Nonetheless, we would like to point out that optimal transport is a useful tool for registration problems, see for example \cite{feydy2017optimal, de2023diffeomorphic,shen2021accurate}.

\subsection{Slice-matching maps}

In this paper, we are interested in maps defined by a slicing and matching procedure \cite{pitie2007automated}, which is closely related to the sliced Wasserstein distance. These are maps of the form
\begin{equation*}
    T_{\sigma,\mu;P}(x) = \sum_{i=1}^n T_{\sigma^{\bth_i}}^{\mu^{\bth_i}}(\bx \cdot \bth_i)\, \bth_i, \quad x\in \mathbb{R}^n
\end{equation*}
involving a source measure $\sigma$, a target measure $\mu$ and an orthogonal matrix $P=[\theta_1,\ldots,\theta_n]$. Here $\sigma^{\bth_i}$ denotes the $1$-dimensional measure obtained by projecting $\sigma$ onto the line $\theta$, and $T_{\sigma^{\bth_i}}^{\mu^{\bth_i}}$ is the optimal transport map between the $1$-dimensional measures $\sigma^{\bth_i}$ and $\mu^{\bth_i}$. Note that the $1$-dimensional optimal transport map can be computed explicitly, see \Cref{sec:prelim} for more details.

This maps allows to define the iterative  \emph{slice-matching procedure} \cite{pitie2007automated,li2023measure}
\begin{equation}\label{intro:iterative_scheme}
    \sigma_{k+1} = (T_{\sigma_k,\mu;P_k})_{\sharp}\sigma_k, \quad k\geq 0,
\end{equation}
which have been successfully used in applications such as color transfer \cite{pitie2007automated}, texture mixing \cite{rabin2012wasserstein} and shape retrieval \cite{rabin2010ShapeRetrieval}. 
Convergence results of \eqref{intro:iterative_scheme} in special cases have been obtained in \cite{pitie2007automated, bonnotte13thesis}. More general almost sure (a.s.) convergence of $\sigma_k \to \mu$ for a stochastic variant of \eqref{intro:iterative_scheme} have been established in \cite{li2023measure}. 

The procedure \eqref{intro:iterative_scheme} also defines an iteration on the level of maps through $T_{\sigma_k,\mu;P_k}, k\geq 0$ (though the mentioned convergence results only hold for measures). Note that bounds which are valid for maps directly carry over to the associated pushforward measures through the well-known stability result \cite[Equation (2.1)]{Ambrosio2013}
\begin{equation}\label{eq:intro-stability}
    W_2(F_{\sharp}\sigma,G_{\sharp}\sigma) \leq \|F-G\|_{\sigma}.
\end{equation}

In this paper, we are interested in the approximation power of one step of \eqref{intro:iterative_scheme}, both for measures and maps. This means that we study (A) the relation between $\sigma_1 = (T_{\sigma,\mu;P})_{\sharp}\sigma$ and the target $\mu$ and (B) the relation between $T_{\sigma,\mu;P}$ and $T$, with $\mu=T_{\sharp}\sigma$.

\subsection{Main contributions}

The contributions of this paper are twofold and summarized in \Cref{intro:one-step-recovery} and \Cref{intro:theorem-affine}. To formulate our contributions, we use the notation $\affset = \{x \mapsto ax+b: a>0, b\in \mathbb{R}^n\}$ for the set of shifts and scaling, and 
$\mathfrak{S}(P)=\{x\mapsto \sum_{i=1}^n f_i(x\cdot\theta_i)\theta_i: f_i:\R\rightarrow\R ~\textrm{increasing}\}$ with $P\in O(n)$ for the set of $P$-compatible maps \cite{khurana2022supervised}. Note that $\affset \subset \mathfrak{S}(P)$ for any $P$.

\begin{theorem}[Recovery and approximation: informal implications of  \Cref{prop:one_step_recovery}, \Cref{prop:Lipbounds}, and \Cref{prop:perturbation}]\label{intro:one-step-recovery}
Consider two measures $\sigma,\mu$ with $\mu=T_{\sharp}\sigma$. Then we get
\begin{enumerate}
\item If $\|T-S\|\leq \varepsilon$ for some $S\in \affset$, then one step of the scheme \eqref{intro:iterative_scheme} reconstructs $T$ up to an error of order $\varepsilon$, i.e. $\|T-T_{\sigma,\mu;P}\|\leq 2\varepsilon$ for any $P \in O(n)$. In particular, if $T\in \affset$, then $T_{\sigma,\mu;P}=T$ for any $P \in O(n)$.
\item If $T\in \mathfrak{S}(P)$ for some orthogonal matrix $P$, then one step of the scheme \eqref{intro:iterative_scheme} using $Q$ reconstructs $T$ up to an error of order $\|P-Q\|_F$, i.e. $\|T-T_{\sigma,\mu;Q}\|\leq C \|P-Q\|_F$. In particular, if $Q=P$, then $T_{\sigma,\mu;Q}=T$.
    \end{enumerate}
\end{theorem}
Through the stability bound of \eqref{eq:intro-stability}, all results of \Cref{intro:one-step-recovery} also hold for the reconstruction of the target measure $\mu$ through $\sigma_1$.

Our results shows that basic transformations relating $\sigma$ to $\mu$ can be recovered easily, not needing any optimization scheme. This relates to recent efforts in trying to approximate the optimal transport map (or fully replace it) by simpler maps such as the Knothe-Rosenblatt construction \cite{baptista2023ApproxFramework}.

To formulate our second contribution, we define the \emph{slice-matching operator} $\smop$, which assigns the first step of \eqref{intro:iterative_scheme} to a given source $\sigma$, target $\mu$, and slicing directions $P\in O(n)$:
\begin{equation*}
    \smop: (\sigma,\mu,P) \mapsto (T_{\sigma,\mu;P})_{\sharp}\sigma.
\end{equation*}

\begin{theorem}[Encoding of special affine effects: informal implications of \Cref{prop:P_inv_equiv}, \Cref{prop:affine_approx}, \Cref{prop:equal_mean}, \Cref{Cor:W2vsSW2}] \label{intro:theorem-affine}
Consider a source measure $\sigma$ and a target measure $\mu$.
One step of the slice-matching procedure \eqref{intro:iterative_scheme} encodes basic transformations in the following sense:
\begin{enumerate}
    \item $\smop$ is invariant to  $\mathfrak{S}(P)$-actions on the source measure $\sigma$:
    \begin{equation*}
    \smop(T_{\sharp}\sigma,\mu,P) = \smop(\sigma,\mu,P), \quad T\in \mathfrak{S}(P).
    \end{equation*}
    \item $\smop$ is equivariant to  $\mathfrak{S}(P)$-actions on the target measure $\mu$:
    \begin{equation*}
    \smop(\sigma,T_{\sharp}\mu,P) = T_{\sharp}\smop(\sigma,\mu,P), \quad T\in \mathfrak{S}(P).
    \end{equation*}
    \item $\smop$ encodes translation effects between $\sigma$ and $\mu$ by matching means: $E(\smop(\sigma,\mu,P))=E(\mu)$ for any $P \in O(n)$.
    \item\label{it:transl-scaling-effect} $\smop$ encodes translation-and-scaling effects in the following sense: Let $S$ be the best map in $\affset$ that aligns $\sigma$ and $\mu$, and let $S^{\ast}$ be the best map in $\affset$ that aligns $\sigma$ and $\smop(\sigma,\mu,P)$. Then, by choosing $P$ randomly, in expectation we get
    \begin{equation*}
       \mathbb{E}\|S-S^{\ast}\|_{\sigma} =C_{\sigma}\left(W_2^2(\sigma,\mu)-nSW^2_2(\sigma,\mu)\right)\geq 0,
    \end{equation*}
    where $C_{\sigma}$ depends on the mean and second moment of $\sigma$, $W_2$ denotes the Wasserstein distance, and $SW_2$ denotes the sliced Wasserstein distance.
\end{enumerate}
\end{theorem}
Note that \eqref{it:transl-scaling-effect} means the following: If $(W_2^2(\sigma,\mu)-nSW^2_2(\sigma,\mu))\leq \varepsilon$, then one step of the iterative scheme \eqref{intro:iterative_scheme} removes global translation-and-scaling effects between the source $\sigma$ and the target $\mu$ up to an error of size $C_{\sigma}\varepsilon$.

\subsection{Structure of the paper}
% \todo[inline]{label sections and use Cref for referencing}
This paper is organized as follows: \Cref{sec:prelim} provides essential background information on optimal transport. \Cref{sec:slice-matching} delves into the details of slice-matching maps, its relations to compatibility, as well as moment-matching properties. In \Cref{sec:slice-operator}, we present invariance, equivariance, and Lipschitz properties associated with the slice-matching operator, which lead to recovery and stability results. In \Cref{sec:registration}, we further explore how the slice-matching procedure handles affine effects by studying basic registration problems. The paper closes with a concluding remark in \Cref{sec:conclusion}.

\section{Preliminaries}\label{sec:prelim}

%\subsection{Wasserstein and sliced Wasserstein distances}\label{sec:prelim-W2}
We use the notation $\mathcal{P}(\R^n)$ and $\mathcal{P}_{ac}(\R^n)$ for the spaces of probability measures on $\R^n$ and absolutely continuous measures with respect to the Lebesgue measure, respectively. We consider the quadratic Wasserstein space, denoted by $W_2(\R^n)$, which includes probability measures $\sigma$ with finite second moments, i.e.\ $\sigma$ satisfying $M_2(\sigma)=\int_{\R^n}\|x\|^2d\sigma(x)<\infty$. In addition, let $\W_{2,ac}(\R^n)=\W_2(\R^n) \cap \mathcal{P}_{ac}(\R^n)$. The mean of a measure $\sigma$ is denoted by $E(\sigma)=\int xd\sigma(x)$.

On $W_2(\R^n)$ we consider the quadratic Wasserstein distance \cite{Villani1}:
\begin{equation*}%\label{wassDist}
 W_2(\sigma,\mu) := \inf_{\pi\in\Gamma(\sigma,\mu)} \left(\int_{\R^{2n}}\|x-y\|^2d\pi(x,y)\right)^{\frac12},   
\end{equation*}
where $\Gamma(\sigma,\mu):=\{\pi \in \mathcal{P}(\R^n\times \R^n): \pi(A\times\R^n)= \sigma(A), \pi(\R^n\times B)=\mu(B), A, B \subseteq \R^n \textrm{ measurable}\}$ represents the set of couplings between $\sigma$ and $\mu$.

When $\sigma \in \W_{2,ac}(\R^n)$ and $\mu\in \W_2(\R^n)$, the optimization problem:
\[
\min_{T:T_{\sharp}\sigma=\mu}\int_{\R^n}\|T(x)-x\|^2 \, d\sigma(x),
\]
with $T$ a map in $L^2(\R^n,\sigma)$, has a unique (up to constants) solution \cite{brenier1991}, which we denote by $T_{\sigma}^{\mu}$. Here $\sharp$ is the pushforward operator. The map $T_{\sigma}^{\mu}$ takes the form $T_{\sigma}^{\mu}=\nabla \varphi$ where $\varphi$ is convex \cite{brenier1991}. Maps which are the gradients of convex functions will be referred to as \textit{Brenier maps}.

If $T_{\sigma}^{\mu}$ exists, the optimal coupling has the form $\pi = (\operatorname{id},T_{\sigma}^{\mu})_{\sharp}\sigma$. In this case, the Wasserstein-2 distance can then be written as:
$
W_2(\sigma,\mu)= \|T_{\sigma}^{\mu}-\operatorname{id}\|_{\sigma},
$
where $\|\cdot \|_{\sigma}$ is the $L^2$-norm w.r.t.\ $\sigma$.

For $1$-dimensional measures, the exist explicit formulae for the optimal transport map and the Wasserstein distance. With $\sigma \in\mathcal{P}_{ac}(\R)$ and $\mu \in \mathcal{P}(\R)$ we get
$
T_{\sigma}^{\mu} =  F_{\mu}^{-1}\circ F_{\sigma},
$
where $F_{\sigma}$ is the cumulative distribution function (CDF) of $\sigma$, and $F_{\mu}^{-1}$ is the pseudo-inverse of the CDF of $\mu$. This leads to:

\begin{equation}\label{eq:one-d-W2}
W_2(\sigma,\mu) =\left( \int_{0}^1 |F_{\mu}^{-1}(x)-F_{\sigma}^{-1}(x)|^2\,dx\right)^{1/2}.
\end{equation}

Throughout this paper, we use the same symbols to denote Wasserstein distance and optimal transport maps for probability measures on both $\R^n$ and $\R$, with the context clarifying the dimension of the measures.

We also use the \emph{sliced-Wasserstein distance} between $\sigma \in \W_{2,ac}(\R^n)$ and $\mu\in \W_2(\R^n)$:
\begin{equation}\label{intro:cont-SW2}
    SW_2^2(\sigma,\mu) = \int_{S^{n-1}} W_2^2(\sigma^{\theta},\mu^{\theta}) \, du(\theta),
\end{equation}
with $\sigma^{\theta} = {\mathcal{P}_{{\theta}}}_{\sharp}\sigma$,
where ${\mathcal{P}_{{\theta}}}(x)= x\cdot \theta$ denotes the projection onto the line defined by $\theta$, and $u$ denotes the uniform measure on $S^{n-1}$. In \eqref{intro:cont-SW2}, $W_2$ denotes the Wasserstein distance between the $1$-dimensional projected measures $\sigma^{\theta},\mu^{\theta}$.

It is known that $SW_2 \leq W_2$ while these two distances are equivalent for measures with compact supports \cite{bonnotte13thesis}.

\section{Slice-matching maps, compatibility and moment matching}\label{sec:slice-matching}
Slice-matching schemes were first introduced by \cite{pitie2007automated} to iteratively transport an initial measure to target measure. An almost sure convergence result of such iterative schemes has been shown in \cite{li2023measure}. In this paper, we are interested in approximation properties of one step of this slice-matching procedure. In what follows, we present the definitions of the schemes and the associated \emph{slice-matching maps}. {We furthermore show that the mean and second moments of $\sigma$ and $\mu$ are matched through one step of such schemes}.

\subsection{Slice-matching maps and compatibility}\label{sec:slice-matching maps}

\begin{definition}[Single-slice and matrix-slice matching, \cite{pitie2007automated, li2023measure}]\label{def:single-matching}
Consider {$\sigma \in \W_{2,ac}(\R^n)$, $\mu\in \W_2(\R^n)$} and a vector $\theta \in S^{n-1}$. The \newmap map from ${\sigma}$ to ${\mu}$ is defined by
\begin{equation}\label{eq:single-Sliceform}
    T_{\sigma,\mu;\theta}(\bx) = \bx + (T_{\sigma^{\theta}}^{\mu^{\theta}}(x\cdot \theta)-x\cdot \theta)\,\theta
\end{equation}
where $T_{\sigma^{\bth}}^{\mu^{\bth}}$ is the optimal transport map between the $1$-dimensional measures $\sigma^{\bth}$ and $\mu^{\bth}$ obtained through projection by $\bth$. If an orthonormal basis of $\R^n$ is used,  the matrix-slice matching map from ${\sigma}$ to ${\mu}$ is defined by
\begin{equation}\label{eq:Sliceform}
    T_{\sigma,\mu;P}(\bx) = \bx + P\begin{bmatrix}T_{\sigma^{\bth_1}}^{\mu^{\bth_1}}(\bx\cdot \bth_1)-\bx\cdot \bth_1\\T_{\sigma^{\bth_2}}^{\mu^{\bth_2}}(\bx\cdot \bth_2)-\bx\cdot \bth_2\\ \vdots\\T_{\sigma^{\bth_n}}^{\mu^{\bth_n}}(\bx\cdot \bth_n)-\bx\cdot \bth_n,
    \end{bmatrix} = \sum_{i=1}^n T_{\sigma^{\bth_i}}^{\mu^{\bth_i}}(\bx \cdot \bth_i)\, \bth_i
\end{equation}
where $P = [\bth_1, \cdots, \bth_n]$ is an orthogonal matrix.
\end{definition}

\begin{remark}\label{rmk:slice-matching property}
The  motivation for the name slice-matching map  is the following: If $\nu = (T_{\sigma,\mu;P})_\sharp \sigma$, then $\nu^{\theta_i}=\mu^{\theta_i}$ for $1\leq i \leq n$, i.e.\ all slices are matched. Similar properties hold for $T_{\sigma,\mu;
\theta}$. Moreover, the following relation between $n$-dimensional Wasserstein distance and one-dimensional Wasserstein distance of the corresponding slices holds:
\begin{equation}\label{eq:WassDist_relation}
    W_2^2(\sigma,(T_{\sigma,\mu;P})_\sharp \sigma) = \sum_{i=1}^nW_2^2(\sigma^{\theta_i},\mu^{\theta_i}),
\end{equation}
see \cite[Lemma 3.9]{li2023measure}. An analogous result  for empirical measures can be found in \cite[ Proposition 5.2.7]{bonnotte13thesis}.
\end{remark}
\begin{remark}\label{remark:iterative-schemes}
The matrix-slice matching maps (as well as the single-slice and generalizations to $1\leq j\leq n$ slices, see \cite{li2023measure}) can be used to approximate a target measure $\mu$ by iteratively pushing a source measure $\sigma_0$:
\begin{equation}\label{eq:iter-scheme}
    \sigma_{k+1} = ((1-\gamma_k)\id + \gamma_k\,{T_{\sigma_k,\mu;P_k}})_{\sharp}\sigma_k, \quad k\geq 0,
\end{equation}
where $\gamma_k$ is a sequence of step-sizes and $P_k$ are matrices in $O(n)$. When $\gamma_k$ satisfies classical stochastic gradient descent assumptions and $P_k$ are chosen i.i.d.\ form the Haar measure on $O(n)$ (and some technical details are satisfied), then $\sigma_k \to \mu$ in both $W_2$ and $SW_2$ a.s. \cite{li2023measure}.

The paper \cite{pitie2007automated} considers the above iterative scheme with $\gamma_k=1$, whose convergence is however not covered by results of \cite{li2023measure}. \cite{pitie2007automated} shows convergence for special measures (the target is Gaussian) and in the KL-divergence.

In this paper, we study approximation properties of one step of \eqref{eq:iter-scheme} with $\gamma_k = 1$, i.e.\ we are interested in the relation between $\sigma_1 = (T_{\sigma_0,\mu;P})_{\sharp}\sigma_0$ and the target $\mu$. An illustration of such approximations using different orthogonal matrices $P$ is given in \Cref{subfig:right}.
%and cases when $W_2(\sigma_1,\mu)$ is small.
\end{remark}

%\begin{remark}\label{rmk:compatible}
 It follows easily that both $T_{\sigma,\mu;\theta}$  and $ T_{\sigma,\mu;P}$ are Brenier maps, though they are not necessarily optimal  transport maps between $\sigma$ and $\mu$. The maps $T_{\sigma,\mu;\theta}$ and $ T_{\sigma,\mu;P}$ are furthermore a  special type of \emph{$P$-compatible maps}  \cite{khurana2022supervised,aldroubi20}, which are defined by
  \begin{equation}\label{eq:compatibleSets}
     \mathfrak{S}(P)=\Big\{x \mapsto \sum_{i=1}^n f_i(x\cdot\theta_i)\theta_i: f_i:\R\rightarrow\R ~\textrm{is increasing}  \Big\},
 \end{equation}
for any fixed $P = [\theta_1,\ldots,\theta_n]  \in O(n)$. An analog to \cref{eq:WassDist_relation} holds for $P$-compatible maps
\begin{equation*}%\label{eq:WassDist_relation-comp}
    W_2^2(\sigma,T_{\sharp}\sigma) = \sum_{i=1}^n W_2^2(\sigma^{\theta_i}, (T_{\sharp}\sigma)^{\theta_i}), \quad T\in \mathfrak{S}(P).
\end{equation*}
Moreover, the slice-matching map can be viewed as the minimizer in $\mathfrak{S}(P)$  associated with the following minimization problem
\begin{align}\label{eq:semi-min-problem}
    T_{\sigma,\mu;P} &= \argmin_{T\in \mathfrak{S}(P)}\sum_{i=1}^nW_2^2((T_{\sharp}\sigma)^{\theta_i},\mu^{\theta_i}).\\ \nonumber
    & \in  \argmin_{T ~\textrm{is Brenier}} \sum_{i=1}^nW_2^2((T_{\sharp}\sigma)^{\theta_i},\mu^{\theta_i}),
\end{align}
where \eqref{eq:semi-min-problem} follows from the fact that for $T\in \mathfrak{S}(P)$, and therefore
\begin{equation*} \sum_{i=1}^nW_2^2((T_{\sharp}\sigma)^{\theta_i},\mu^{\theta_i}) = \|T-T_{\sigma,\mu;P}\|^2_{\sigma}, ~\textrm{which is minimal iff~} T = T_{\sigma,\mu;P}.
\end{equation*}
The details of this statement are presented in  \Cref{cor:min-set-compatible}.
%\end{remark}

$P$-compatible maps have been used in tangent space embeddings, which allow for linear separability of two classes of measures, see \cite{khurana2022supervised,aldroubi20,park18,moosmueller2020linear}. These maps satisfy
\begin{equation}\label{SW-compa}
    \mathcal{P}_{\theta_i}\circ T  = f_i \circ \mathcal{P}_{\theta_i}, \quad i=1,\ldots, n,
\end{equation}
  where $T(x) = \sum_{i=1}^n f_i(x\cdot\theta_i)\theta_i \in \mathfrak{S}(P)$ and $P=[\theta_1,\ldots,\theta_n]$.

Next, we show that the property \eqref{SW-compa} also characterizes the set of $P$-compatible maps given in \eqref{eq:compatibleSets}.

 \begin{proposition}\label{lem:commute_proj}
Let $T \in \mathfrak{S}(P)$ with $T(x) = \sum_{i=1}^n f_i(x\cdot\theta_i)\theta_i$, where $[\theta_1,\ldots,\theta_n] = P \in O(n)$ and $f_i, i=1,\ldots,n$ are increasing functions. Then $T=T_{\sigma}^{\mu}$ for some measures $\sigma \in \W_{2,ac}(\R^n)$ and $\mu \in \mathcal{P}_{ac}(\R^n)$, if and only if
\begin{equation}\label{eq:commute_proj}
    \mathcal{P}_{\theta_i} \circ T_{\sigma}^{\mu} = f_i \circ \mathcal{P}_{\theta_i}, \quad i=1,\ldots,n.
\end{equation}
Furthermore, in this case, $f_i = T_{\sigma^{\theta_i}}^{\mu^{\theta_i}}$ and $T$ is a \Mnewmap map.
\end{proposition}
\begin{proof}
    For the equivalence, note that relation \eqref{eq:commute_proj} is equivalent to
    \begin{equation*}
        T(x) = \sum_{i=1}^n f_i(\bx \cdot \bth_i)\, \bth_i = \sum_{i=1}^n (\bth_i \cdot T_{\sigma}^{\mu}(\bx))\, \bth_i = T_{\sigma}^{\mu}(\bx),
    \end{equation*}
    since the columns of $P$ are orthonormal.

 For the other statement note that \eqref{SW-compa}
 implies that
\begin{align*}
    \mu^{\theta_i} &= \left(\mathcal{P}_{\theta_i}\circ T\right)_{\sharp} \sigma 
    = \left(f_i \circ \mathcal{P}_{\theta_i}\right)_{\sharp} \sigma \\
    & = (f_i)_{\sharp}\sigma^{\theta_i}.
\end{align*}
Since $f_i$ are increasing we obtain $f_i = T_{\sigma^{\theta_i}}^{\mu^{\theta_i}}$.
\end{proof}

\begin{remark}
    Note that the proof of \Cref{lem:commute_proj} relies on the columns of $P$ forming an orthonormal basis of $\R^n$. Therefore, this equivalence is not true for \newmap maps. We only have the implication $T_{\sigma,\mu;\theta} = T_{\sigma}^{\mu} \implies \mathcal{P}_{\theta} \circ T_{\sigma}^{\mu} = T_{\sigma^{\theta}}^{\mu^{\theta}} \circ \mathcal{P}_{\theta}$.
\end{remark}

\subsection{Moment matching}
We show that slice-matching maps push the source measure to a measure that has the same mean and second moments as the target measure. 
\begin{proposition}
Let $\sigma\in \W_{2,ac}(\R^n)$ and $\mu\in \W_2(\R^n)$. Then for any $P\in O(n)$, the following holds: 
   \begin{equation}\label{eq:equalmean}
      E((T_{\sigma,\mu;P})_\sharp \sigma) = E(\mu)
   \end{equation} 
   \begin{equation}\label{eq:equal2ndMoment}
       M_2((T_{\sigma,\mu;P})_\sharp \sigma) = M_2(\mu)
   \end{equation}
\end{proposition}

\begin{proof}
A direct computation shows that  $ \int T_{\sigma,\mu;P}(x) d\sigma(x) = \int yd\mu(y)$ and $\int \|T_{\sigma,\mu;P}(x)\|^2 d\sigma(x) = \int \|y\|^2d\mu(y)$, see \Cref{lm:equal_mean} for more details.
\end{proof}

The equal mean property \eqref{eq:equalmean} gives a first hint towards the shift-eliminating phenomenon of the slice-matching procedure. In particular, it can be verified that if $T_{\sigma}^{\mu}$ is a shift, then $(T_{\sigma,\mu;P})_\sharp \sigma = \mu$.
 A more comprehensive perspective of the shift-eliminating effect will be presented in  \cref{S_equiVariance} and \cref{equal-shift-registration}.

\section{Invariance, equivariance, and Lipschitz properties}\label{sec:slice-operator}
We consider the following operator induced by the slice-matching maps of \Cref{def:single-matching}:

\begin{definition}[Slice-matching operator]\label{def:smop}
We define the following operator based on the slice-matching approximation,  concerning a source measure, a target measure, and slicing directions given by an orthogonal matrix:
\begin{align*}
       \smop:\ \W_{2,ac}(\R^n)\times \W_{2}(\R^n)\times O(n)&\rightarrow \W_{2}(\R^n)\\
       \quad \quad    (\sigma,\ \mu,\ P)&\mapsto (T_{\sigma,\mu;P})_{\sharp}\sigma.
   \end{align*}
\end{definition}
Note that $(T_{\sigma,\mu;P})_{\sharp}\sigma$ has finite second moment when $\mu\in \W_2(\R^n)$ by \eqref{eq:equal2ndMoment}, which implies that $\smop$ maps into $\W_{2}(\R^n)$.
If absolute continuity of $(T_{\sigma,\mu;P})_{\sharp}\sigma$ is  desired, one can further assume that both $\sigma,\mu$ are absolutely continuous, in which case $\smop: \W_{2,ac}(\R^n)\times \W_{2,ac}(\R^n)\times O(n) \rightarrow \W_{2,ac}(\R^n)$, see \Cref{lm:abs_cont}.

We first illustrate the invariance and equivariance properties of the slice-matching operator in terms of basic transformations, i.e., shifts and scalings, which can be viewed as special cases of compatible transformations, as shown in \Cref{sec:CompInvar}. {We also show how such properties are related to the recovery of optimal transport maps using matrix-slice matching. Moreover, as a complementary remark to the unifying convergence framework \cite{li2023measure} of the single-slice and matrix-slicing matching schemes,  we illustrate their differences via different recovery properties.}  A Lipschitz property in terms of the third component is shown in \Cref{sec:Lip}.

\subsection{Invariance and equivariance with respect to shifts and scalings}
We show that the {\opname}  is invariant to actions induced by push-forward operations of shifts and scalings on any initial measure $\sigma\in \W_{2,ac}(\R^n)$ and is equivariant to actions of these maps on any target measure $\mu \in \W_2(\R^n)$, regardless of the orthogonal matrix $P$. More specifically, 
\begin{proposition}\label{prop:P_inv_equiv}
Let $\sigma\in \W_{2,ac}(\R^n)$ and $\mu\in \W_2(\R^n)$. Then for any $P\in O(n)$ and $S(x)=ax+b$ where $a>0, b\in \R^n$,  we get
   \begin{align}
        \smop(S_{\sharp}\sigma,\mu,P) &= \smop(\sigma,\mu, P),\label{S_invariance} \\
          \smop(\sigma,S_\sharp\mu,P) &= S_\sharp\smop(\sigma,\mu, P).\label{S_equiVariance}
   \end{align}
\end{proposition}

\begin{proof}
   Since $S\in \bigcap\limits_{P\in O(n)}\mathfrak{S}(P)$, the conclusion follows from the corresponding invariance and equivariance properties in terms the group $\mathfrak{S}(P)$ of  compatible transformations, see \Cref{prop:invariant_Porbits}.
\end{proof}

\begin{figure}
    \centering

    \begin{subfigure}[b]{0.48\textwidth}
        \centering
        \includegraphics[width=1.05\textwidth]{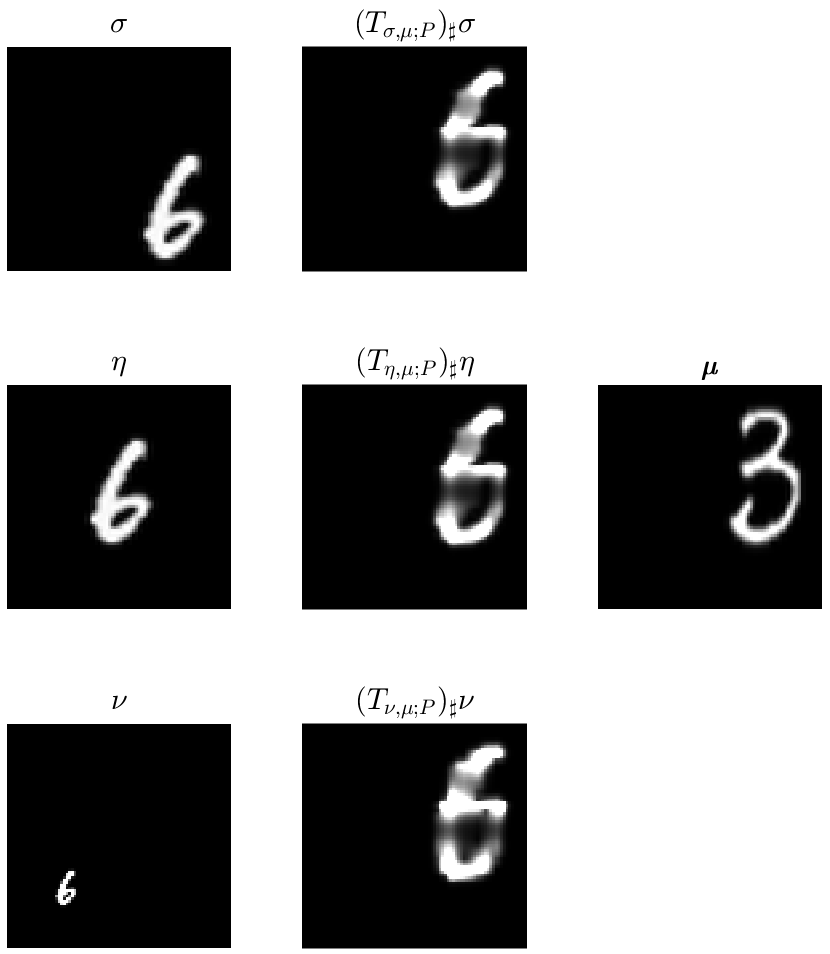}
        \caption{First column: different initial images $\sigma, \eta, \nu$ created by translating and scaling a particular image of the digit $6$. Middle column: images $(T_{\sigma_i,\mu;P})_{\sharp}\sigma_i$ generated via slice-matching, where $P= \begin{bmatrix}
            1&0\\0&1
        \end{bmatrix}.$ Third column: a fixed target image of the digit 3.}
        \label{subfig:left}
    \end{subfigure}
\hfill
    % Subfigure 2 (right)
    \begin{subfigure}[b]{0.44\textwidth}
        \centering
        \includegraphics[width=1.1\textwidth]{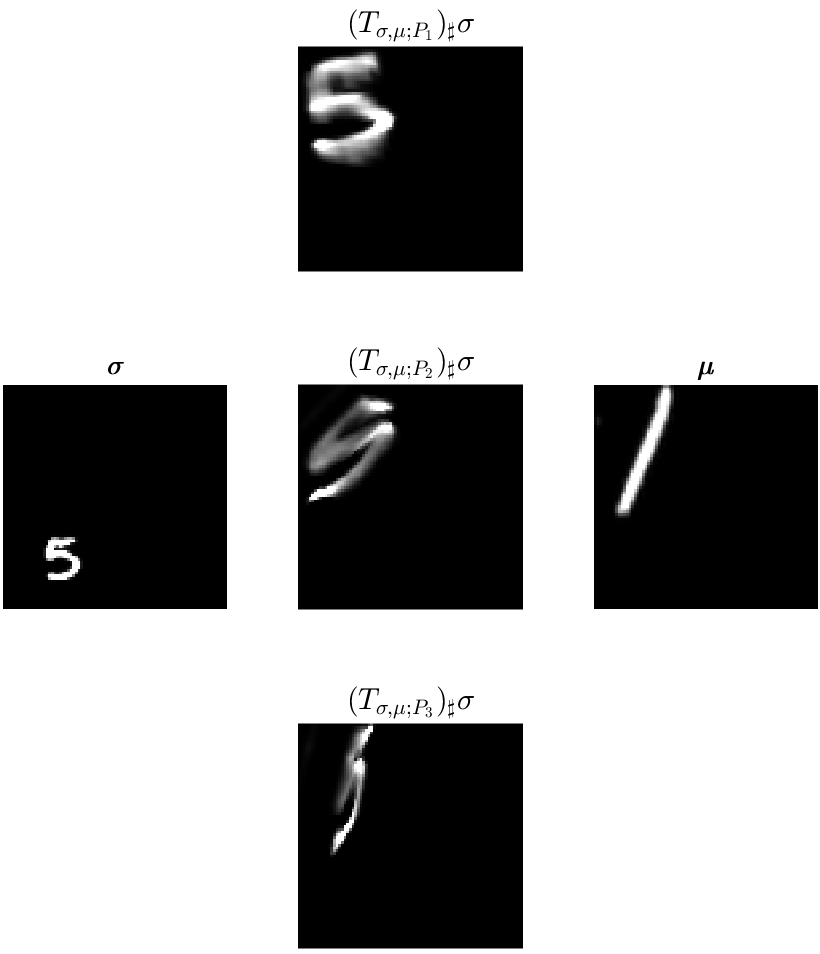}
        \caption{First and third columns: initial image $\sigma$ of the digit 5 and target image $\mu$ of the digit 1 . Middle column: images $(T_{\sigma,\mu;P_i})_{\sharp}\sigma$ generated via slice-matching with different orthogonal matrices $P_i = \begin{bmatrix}
            \cos\theta_i& \sin\theta_i\\-\sin\theta_i&\cos\theta_i
        \end{bmatrix}$, where $\theta_i = \frac{(2i-1)\pi}{12}, i = 1,2,3$.}
        \label{subfig:right}
    \end{subfigure}

    \caption{Effects of matrix-slice matching. Left panel: Illustration of invariance to translation and scaling of the initial measures. Right panel: Illustration of slice-matching using different orthogonal matrices.}
    \label{fig:SliceMatching-effects}
\end{figure}
An illustration of the invariance property the $\smop$ with respect to translation-and-scaling transformations of the source measure is presented in \Cref{subfig:left}. Approximations with different choices of orthogonal matrices is illustrated in \Cref{subfig:right}.

Note that isotropic scalings and translations $S(x)= a\bx +\bb$ with $a>0$ and $b \in \mathbb{R}^n$ are special types of compatible maps. They satisfy $S\in \mathfrak{S}(P)$ for all $P \in O(n)$. An immediate corollary of the above proof is that, given two measures that are related by shifts and scalings, the target is recovered exactly by push the initial measure with the slice matching map $T_{\sigma,\mu;P}$ for any $P$.  
\begin{corollary}[Recovery of basic transformations with one step of slice matching map]\label{prop:one_step_recovery}
Given $\sigma\in \W_{2,ac}(\R^n),\mu\in \W_2(\R^n)$ with $T_{\sigma}^{\mu}(x)=ax+b$ for some $a>0,b\in\R^n$. Then we have $T_{\sigma,\mu;P} = T_{\sigma}^{\mu}$ and $\smop(\sigma,\mu,P)=\mu$ for any $P\in O(n)$.
\end{corollary}
\begin{proof}
 Let $S=T_{\sigma}^{\mu}$ in \eqref{S_invariance}. Then   $S_{\sharp}\sigma = \mu$ and $\smop(S_{\sharp}\sigma,\mu,P) = (T_{\mu,\mu;P})_{\sharp}\mu = (\operatorname{id})_{\sharp}\mu=\mu$. Hence by  \eqref{S_invariance}, $\smop(\sigma,\mu, P)=\mu$. The fact that $T_{\sigma,\mu;P}= T_{\sigma}^{\mu}$ follows from the fact that they are both Brenier maps pushing $\sigma$ to $\mu$. 
\end{proof}

\begin{figure}
    \centering
    \includegraphics[width=.5\textwidth]{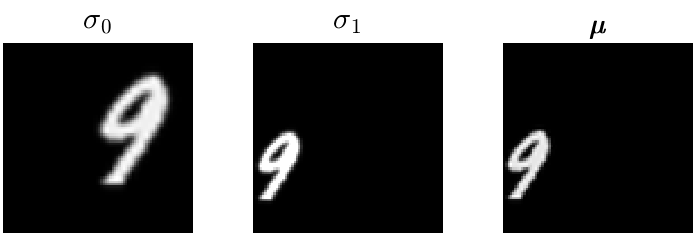}
    \caption{Recovery of basic transformation (shift and scaling) using matrix-slice matching. The initial $\sigma_0$ (left) and target image $\mu$ (right)   are related $T_{\sharp}\mu=\sigma_0$ where $T(x)= 1.6(x+[-35,20]^t)$, where each image is of size $84\times 84$. Here $\sigma_1= (T_{\sigma_0,\mu;P})_\sharp\sigma_0$ (middle) is the image obtained via slice-matching transport between $\sigma_0$ and $\mu$ for an arbitrary orthogonal matrix $P$.  Note $\sigma_1\approx \mu$ up to numerical errors, as implied by \Cref{prop:one_step_recovery}.}
    \label{fig:Exp2}
\end{figure}

An illustration of \Cref{prop:one_step_recovery} is presented in \Cref{fig:Exp2}, the target image $\mu$ is matched (up to numerical errors) by its slice-matching approximation $(T_{\sigma,\mu;P})_\sharp\sigma$ for any $P$, if $T_{\sigma}^{\mu}$ is a translation-scaling function.

  In addition, we can show that a differentiable map $T$ connecting $\sigma$ and $\mu$ can only be recovered with one step of the slice-matching scheme with any choice of $P$ if and only if $T$ is an isotropic scaling with translation:

  \begin{remark}\label{recovery_iff}
Under the assumptions in \Cref{prop:one_step_recovery}, and if we further assume that $T_{\sigma}^{\mu}$ is differentiable, we obtain the following: $T = T_{\sigma,\mu;P}$ for any choice of $P \in O(n)$ if and only if $T(x)=ax+b$ for some $a>0$ and $b\in \R^n$. One direction follows from \Cref{prop:one_step_recovery}. See \Cref{prop:intersect_compatible} for details of the other direction.  
  \end{remark}

The above recovery result holds for the  matrix-slice scheme, but in general does not hold for single-slice schemes:

\begin{example}\label{example_singleSlice}
Let  $\mu = T^b_{\sharp}\sigma$ with $T^b(x)= x+b, b\neq 0\in \R^n$. Unlike the \Pmatching, the target measure $\mu$ cannot always be recovered via $(T_{\sigma,\mu;\theta})_\sharp\sigma$ by the observing that
 $T_{\sigma,\mu;\theta}(x) = x+ \theta(\theta\cdot b),$ for $\theta\in S^{n-1}$. To recover $\mu$ exactly, $\theta=\frac{1}{\|b\|}b$ is the only choice. Therefore, in general, we will not recover $\mu$ after one step.
\end{example}
\begin{discussion}
A possible remedy is through an iterative scheme by repeating the above slice-matching procedure with different $\theta$ as introduced in \Cref{remark:iterative-schemes}. In \cite{li2023measure} we show that the iterative scheme $\sigma_{k+1}=((1-\gamma_k)\operatorname{id}+\gamma_k T_{\sigma_k,\mu; \theta_k})_\sharp \sigma_k$ with step size $\gamma_k$ satisfying $\sum\gamma_k = \infty$ and $\sum \gamma_k^2<\infty$ converges a.s. ($\theta_k\overset{\textrm{i.i.d}}{\sim} u $, $u$ uniform measure on $S^{n-1}$).

However, for the particular example of a shift, a.s.\ convergence can also be achieved for $\gamma_k=1$. In \Cref{proofTranslationConv} we give an elementary proof that $\sigma_{k+1}=(T_{\sigma_k,\mu; \theta_k})_\sharp \sigma_k$ converges to $\mu=T^{b}_{\sharp}\sigma$ a.s. with respect to the $W_2$-distance.
\end{discussion}

\subsection{Invariance and equivariance with compatible maps}\label{sec:CompInvar}
In this subsection, we discuss the invariance and equivariance properties of the slice-matching operator $\smop$ as defined in \Cref{def:smop} concerning compatible transformations, defined in \eqref{eq:compatibleSets}.
\begin{proposition}\label{prop:invariant_Porbits}
Let $\sigma\in \W_{2,ac}(\R^n)$ and $\mu\in \W_2(\R^n)$.
For any $T\in \mathfrak{S}(P)$, where $P\in O(n)$, we have
\begin{align}
      \smop({T}_{\sharp}\sigma,\mu, P)& = \smop(\sigma,\mu, P),\label{eq:P-inv}\\
        \smop(\sigma,{T}_{\sharp}\mu, P)& ={T}_{\sharp} \smop(\sigma,\mu, P). \label{eq:P-eqinv}  
\end{align}
\end{proposition}
\begin{proof}
    Let $T(x) = Pf(P^tx)$, where $f(x)= (f_1(x_1),\cdots, f_n(x_n))$ with each $f_i$ being increasing and $P = [\theta_1,\cdots, \theta_n]$. Let $\sigma_T = T_{\sharp}\sigma$.
    By \eqref{SW-compa}, we get $\sigma_T^{\theta_i}= {f_i}_{\sharp}\sigma^{\theta_i}$. It follows from the fact that $f_i$ is increasing that 
$T_{\sigma^{\theta_i}}^{\mu^{\theta_i}}=T_{\sigma_T^{\theta_i}}^{\mu^{\theta_i}}\circ f_i$. Denote $g(x)= [(T_{\sigma_T^{\theta_1}}^{\mu^{\theta_1}}\circ f_1)(x_1), \cdots, (T_{\sigma_T^{\theta_n}}^{\mu^{\theta_n}}\circ f_n)(x_n)]^t$ and  $g_T(x)= [T_{\sigma_T^{\theta_1}}^{\mu^{\theta_1}}(x_1), \cdots, T_{\sigma_T^{\theta_n}}^{\mu^{\theta_n}}(x_n)]^t$. Then
\begin{align*}
    (T_{\sigma,\mu;P})_{\sharp}\sigma &= (Pg\circ P^t)_{\sharp}\sigma=(Pg_T\circ f\circ P^t)_{\sharp}\sigma\\
    &=(pg_T\circ P^t\circ P\circ f\circ P^t)_\sharp \sigma =(T_{\sigma_T,\mu;P}\circ T)_{\sharp}\sigma\\
    &= (T_{\sigma_T,\mu;P})_{\sharp}\sigma_T.
\end{align*}
This proves \eqref{eq:P-inv}. For \eqref{eq:P-eqinv}, denote $\mu_T = T_{\sharp}\mu$. With similar reasoning as before and with the observations that $T_{\sigma^{\theta_i}}^{\mu_T^{\theta_i}}=f_i\circ T_{\sigma^{\theta_i}}^{\mu^{\theta_i}}$, we get
   \begin{align*}
       T_{\sigma,\mu_T;P} &= Pf\circ P^t\circ T_{\sigma,\mu;P} = T\circ T_{\sigma,\mu;P}. \qedhere
   \end{align*}
\end{proof}
%\mathfrak{M}^P_{\eta}
As a direct conclusion, define $\mathfrak{S}(P)_{\sharp}\eta:= \{T_{\sharp}\eta: T\in \mathfrak{S}(P) \}$, then $\smop(\mathfrak{S}(P)_{\sharp}\sigma,\mu,P)=\{\smop(\sigma, \mu, P)\}$ is a singleton set, as illustrated in \Cref{subfig:left}, and $\smop(\sigma,\mathfrak{S}(P)_{\sharp}\mu,P) = \mathfrak{S}(P)_{\sharp}{\smop(\sigma,\mu,P)}$.

%%%
\begin{corollary}\label{cor:min-set-compatible}
Let $\sigma,\mu\in \W_{2,ac}(\R^n)$. Then
\begin{equation*}   
\sum_{i=1}^nW_2^2((T_{\sharp}\sigma)^{\theta_i},\mu^{\theta_i}) = \|T-T_{\sigma,\mu;P}\|^2_{\sigma},
\end{equation*}
where $T\in \mathfrak{S}(P)$ and $P\in O(n)$.
\begin{proof}
  Rewriting \eqref{eq:WassDist_relation}, we have
  \begin{equation*}
   \sum_{i=1}^n W_2^2(\sigma^{\theta_i},\mu^{\theta_i})= W_2^2(\sigma,\smop(\sigma,\mu,P)).
\end{equation*} 
Hence
\begin{align*}
      \sum_{i=1}^n W_2^2((T_\sharp\sigma)^{\theta_i},\mu^{\theta_i}) &= W_2^2(T_{\sharp}\sigma,\smop(T_{\sharp}\sigma,\mu,P))
       = W_2^2(T_{\sharp}\sigma,\smop(\sigma,\mu,P))\\
      & = W_2^2(T_{\sharp}\sigma, (T_{\sigma,\mu;P})_{\sharp}\sigma)
      = \|T-T_{\sigma,\mu;P}\|^2_{\sigma},     
\end{align*}
where the second step uses the invariance property \eqref{eq:P-inv} and the last steps makes uses of the isometry property with respect to $P$-compatible maps $T$ and $T_{\sigma,\mu;P}$, see \cite{khurana2022supervised}. 
\end{proof}

\end{corollary}

\subsection{Recovery and stability properties of matrix-slicing matching}\label{sec:Lip}
When two measures are related by a $P$-compatible map, the optimal map between them can be recovered exactly by the ($P$)-matrix-slice matching procedure. This is in contrast to shifts and scalings, where an arbitrary orthogonal matrix can be used for recovery. The following corollary to \Cref{prop:invariant_Porbits} summarized this result.

\begin{corollary}[Recovery of $P$-compatible transformations with one step of $P$-slice-matching]\label{prop:one_step_recovery_compatible}
Given $\sigma\in \W_{2,ac}(\R^n),\mu\in \W_2(\R^n)$ with $T_{\sigma}^{\mu}\in \mathfrak{S}(P)$ for some $P\in O(n)$. Then we have $\smop(\sigma,\mu,P)=\mu$ and $T_{\sigma,\mu;P} = T_{\sigma}^{\mu}$. 
\end{corollary}
\begin{proof}
 With \Cref{prop:invariant_Porbits}, the above result follows from similar arguments as in \Cref{prop:one_step_recovery}.
\end{proof}
To recover a compatible map $T$ with one step of the iteration, \Cref{prop:one_step_recovery} implies that we need to know the orthogonal matrix $P$. The following Lipschitz continuity of the slice-matching operator $\smop$ with respect to $P$ establishes a stability result on the choice of $P$:
\begin{proposition}\label{prop:Lipbounds}
Let $\sigma\in \W_{2,ac}(\R^n)$ and $\mu\in \W_2(\R^n)$.  Assume that there exists $L>0$ such that $T_{\sigma^{\theta}}^{\mu^{\theta}}$ is $L$-Lipschitz on $\mathbb{R}$ for all $\theta\in S^{n-1}$. Then 
   \begin{align}\nonumber
     W_2( \smop(\sigma,\mu, P), \smop(\sigma,\mu,Q)) & = \|T_{\sigma,\mu;P}-T_{\sigma,\mu;Q}\|_{\sigma} \\ \label{eq:Lipbounds}
    & \leq (3L+1)C\|P-Q\|_F,
   \end{align}
   where $C = \max\{M_2(\sigma), M_2(\mu)\}$ and $\|\cdot\|_F$ denotes the Frobenius norm.
\end{proposition}
\begin{proof}
See \Cref{append:proof_Lipbounds}.   
\end{proof}
\begin{remark}
Inequality \eqref{eq:Lipbounds} can be viewed a stability result for one step of the iterative schemes described in \Cref{remark:iterative-schemes}.
 If $T_{\sigma}^{\mu}\in \mathfrak{S}(P)$ for some $P\in O(n)$, then the push-forward measure of $\sigma$ using a slice-matching map associated with $Q\in O(n)$ is within $(3L+1)C\|P-Q\|_F$ in Wasserstein distance to the target $\mu$.  Picking a $Q$ close to $P$ is good enough to obtain an approximation of $\mu$ by $\sigma_1 :=(T_{\sigma,\mu;Q})_\sharp\sigma$. Note that if $Q=P$, then $\sigma_1=\mu$, which also follows from \Cref{prop:one_step_recovery_compatible}.
 %$(T_{\sigma,\mu;Q})_\sharp\sigma$
\end{remark}

The stability result in \Cref{prop:Lipbounds} shows how well the target $\mu$ can be approximated by the slice-matching approximation $(T_{\sigma,\mu;P})_{\sharp}\sigma$ when $\sigma$ and $\mu$ are related by some compatible map $T_{\sigma,\mu;Q}$. Additionally, we will now show that if $\sigma$ and $\mu$ are related by a map which is an $\varepsilon$-perturbation of shifts and scalings, then $\mu$ can be approximated by its slice-matching approximation with at most $2\varepsilon$ error (\Cref{rmk:epsilonbound}). This can also be viewed as an extension to the recovery result in \Cref{prop:one_step_recovery}.
\begin{proposition}\label{prop:perturbation}
    Let $\sigma\in \W_{2,ac}(\R^n)$ and $\mu\in \W_2(\R^n)$. Then   for any $P\in O(n)$,
    \begin{equation*}
        W_2((T_{\sigma,\mu;P})_{\sharp}\sigma,\mu)\leq 2 \inf_{S\in \affset}W_2(S_{\sharp}\sigma,\mu),
    \end{equation*}
where $\affset := \{S(x)=ax+b \mid a>0, b\in \R^n\}$.
\end{proposition}
\begin{proof}
Since $W_2(S_\sharp\sigma,\mu)= \|S-T_{\sigma}^\mu\|_{\sigma}$, it suffices to show that for any $S\in \affset$,
\begin{equation*}
     W_2((T_{\sigma,\mu;P})_{\sharp}\sigma,\mu)\leq 2\|S-T_\sigma^\mu\|_{\sigma}.
\end{equation*}
        By the Lipschitz property (see e.g., \cite[Equation (2.1)]{Ambrosio2013}) associated with $W_2$  and triangle inequality, we have
    \begin{align}
W_2((T_{\sigma,\mu;P})_{\sharp}\sigma,\mu) &\leq \|T_{\sigma,\mu;P}- T_\sigma^\mu\|_{\sigma}\nonumber\\
&\leq \|T_{\sigma,\mu;P}- S\|_{\sigma} +\|S- T_\sigma^\mu\|_{\sigma}.\label{ineq-2epsilon}
    \end{align} Next we bound the first term.
Since $S\in \mathfrak{S}(P)$ for any $P\in O(n)$, it follows from  \Cref{cor:min-set-compatible} that 
\begin{align*}%\label{temp-ineq}
    \|S- T_{\sigma,\mu;P}\|^2_{\sigma} &= \sum_{i=1}^n W_2^2((S_{\sharp}\sigma)^{\theta_i},\mu^{\theta_i}) \leq W_2^2(S_{\sharp}\sigma,\mu) = \|S- T_\sigma^\mu\|^2_{\sigma}
\end{align*}
where the bound follows from \Cref{lm:WassSlice_compare} and the last equality follows from isometry properties with respect to transformations in $\affset$ \cite{moosmueller2020linear}. The desired inequality hence follows from \eqref{ineq-2epsilon}.
\end{proof}

\begin{remark}\label{rmk:epsilonbound}
    Assume that  $\mu=(f\circ S)_{\sharp}\sigma $, where $f:\R^n\rightarrow\R^n$ satisfies  $\|f-\id\|_{S_\sharp\sigma}\leq \varepsilon$
    for some $\varepsilon>0$ and $S(x)=ax+b, a>0, b\in \R^n$.
 Then \begin{equation*}
        W_2((T_{\sigma,\mu;P})_{\sharp}\sigma,\mu)\leq 2\varepsilon.
    \end{equation*}
\end{remark}

\begin{remark}
    Using essentially the same arguments as in the proof of \Cref{prop:perturbation}, one can show that for any $T:\R^n\rightarrow \R^n$ such that $\mu= T_\sharp\sigma$,
    \begin{equation*}
      \|T- T_{\sigma,\mu;P}\|_{\sigma} \leq 2\inf_{S\in \affset}\|T-S\|.
    \end{equation*}
Note $T$ is not necessarily the optimal transport map between $\sigma$ and $\mu$.
\end{remark}

\section{Affine effects and registration problems} \label{sec:registration}

We  study two basic image and point-cloud registration problems to understand the effects of the slice-matching maps \eqref{eq:Sliceform}. Image registration \cite{Brown92-survey} involves matching images with variations caused by differences in acquisition, object growth or other changes. It plays a fundamental role in image processing, particularly in medical image applications \cite{Sotiras13survey}. 
Modeling image data and shape with probability measures has paved the way for robust and scalable algorithms by leveraging  the theory optimal transport, such as diffeomorphic registration methods \cite{feydy2017optimal, de2023diffeomorphic}, including point cloud registration \cite{shen2021accurate}.

\begin{figure}
    \centering
    \includegraphics[width=.5\textwidth]{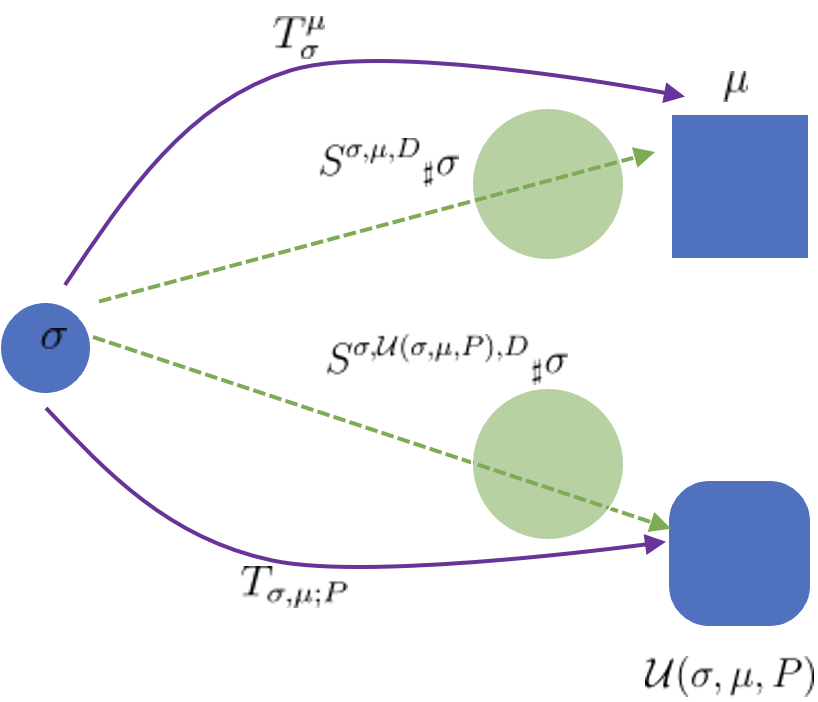}
    \caption{Illustration of registration problems with translation-and-scaling function $S^{\sigma,\eta,D}$ from \Cref{prop:affine_approx}, where $\eta= \mu \textrm{~or~} \eta=\smop(\sigma,\mu,P)$. As indicated in \eqref{eq:registr-equal-mean}, ${S^{\sigma,\eta,D}}_\sharp\sigma$ and $\eta$ have the same mean; however for better visualization, we intentionally kept them separate in this figure.}
    \label{fig:sketch-registration}
\end{figure}

We have shown that slice-matching maps can be used to register translation-and-scaling deformations exactly, see \Cref{prop:one_step_recovery} and \Cref{fig:Exp2}. We also showed that if the two measures are related by perturbations of translations and scalings, the registration error is bounded by the the size of this perturbation, see \Cref{prop:perturbation}. To gain a better understanding of how the slice-matching procedure incorporates affine effects,  particularly when the initial and target measures significantly differ--meaning they are not merely small perturbations of translations and scalings-- we compare the registration maps aimed at the target $\mu$
 with those directed at its slice-matching approximation $\smop(\sigma,\mu,P)$, see \Cref{fig:sketch-registration}. Specifically, we demonstrate that registration maps involving only translations are identical (\Cref{prop:equal_mean}), and those involving translations and isotropic scalings are comparable (\Cref{prop:affine_approx}).

Let $D(\cdot, \cdot)$ be a distance between probability measures, e.g., $W_2$ or $SW_2$. We study  registration problems with the following subsets of affine transformation: $\affset_t := \{S(x)=x+b\mid b\in\R^n\}$ (the set of translations), $\affset := \{S(x)=ax+b \mid a>0, b\in \R^n\}$ (the set of compositions of isotropic scalings and translations):  
\begin{align}
 S_t^{\sigma,\eta,D} &:= \argmin_{S\in \affset_t} D(S_{\sharp}\sigma,\eta),\label{eq:St_registr}\\
    S^{\sigma,\eta,D} &:= \argmin_{S\in \affset} D(S_{\sharp}\sigma,\eta),\label{eq:S_registr}
\end{align}
the existence and uniqueness of which will be addressed later.

We show that the optimal translation registration maps from the initial measure to the target and to slice-matching approximation of the target are identical.

\begin{proposition}[Registration with translations \eqref{eq:St_registr}]\label{prop:equal_mean}
Let $\sigma\in \W_{2,ac}(\R^n)$ and $\mu\in \W_2(\R^n)$. Then for any $P\in O(n)$, the unique minimizer in  \eqref{eq:St_registr} satisfies
\begin{equation}\label{equal-shift-registration}
S_t^{\sigma,\mu,W_2} = S_t^{\sigma,\smop(\sigma,\mu,P),W_2}.
\end{equation}
\end{proposition}
\begin{proof}
The minimization problem \eqref{eq:St_registr} is a quadratic problem in the parameter $b$, which can be solved by taking partial derivatives with respect to $b$, and setting them to $0$. From this, the existence and uniqueness follows immediately. Calculations are summarized in the proof of \Cref{keyproof:prop:affine_approx} and \Cref{keycor:prop:affine_approx}. In particular, the arguments show that the optimal parameters $b^{W_2}, \widetilde{b}^{W_2}$ for $S_t^{\sigma,\mu,W_2}$ and $ S_t^{\sigma,\smop(\sigma,\mu,P),W_2}$ respectively, satisfy
    \begin{equation*}
        b^{W_2}= E(\mu) - E(\sigma), \quad  \widetilde{b}^{W_2} = E(\smop(\sigma,\mu,P)) - E(\sigma).
    \end{equation*}
The conclusion hence follows from the fact that $E(\smop(\sigma,\mu,P)) = E(\mu)$, see \eqref{eq:equalmean}.
\end{proof}
\begin{remark}
    By similar calculations, one can also show that 
    \begin{equation*}
        S_t^{\sigma,\mu,SW_2} = S_t^{\sigma,\smop(\sigma,\mu,P),SW_2}.
    \end{equation*}
    See expressions $b^{SW_2}$ and $ \widetilde{b}^{SW_2}$ in \Cref{keyproof:prop:affine_approx} and \Cref{keycor:prop:affine_approx}.
\end{remark}

The above registration result illustrates the idea of how shifts are eliminated by the slice-matching procedure.  When considering registration involving translations and isotropic scalings measured by the $W_2$ distance, the following comparison holds:
\begin{proposition}[Registration with translation-and-scalings \eqref{eq:S_registr}]\label{prop:affine_approx}
Let $\sigma\in \W_{2,ac}(\R^n)$ and $\mu\in \W_2(\R^n)$. Assume that (i) The convex potential $\phi$ such that $\grad \phi = T_{\sigma}^{\mu}$ given by Brenier's theorem is differentiable at $E(\sigma)$,  and (ii) For any $
    \lambda \in (0,1)$, $\phi((1-\lambda)y+\lambda E(\sigma)) < (1-\lambda)\phi(y)+\lambda\phi(E(\sigma))$ for all $y$ in some  ball $B(x,r)$, where $x$ lies in the support of $\sigma$.  Then $S^{\sigma,\mu,W_2}$ and $S^{\sigma,\smop(\sigma,\mu,P), W_2}$ in \eqref{eq:S_registr} are well-defined and unique, and satisfy the following
    \begin{align}\label{eq:S-Sp}
    W_2({S^{\sigma,\mu,W_2}}_{\sharp}\sigma,{S^{\sigma,\smop(\sigma,\mu,P), W_2}}_{\sharp}\sigma) &=\|S^{\sigma,\mu,W_2} - S^{\sigma,\smop(\sigma,\mu,P), W_2}\|_{\sigma}\nonumber\\& =  \frac{W_2^2(\sigma,\mu) -\sum_{i=1}^nW_2^2(\sigma^{\theta_i},\mu^{\theta_i})}{2\sqrt{M_2(\sigma)-\|E(\sigma)\|^2}},
\end{align}
where $P=[\theta_1,\ldots,\theta_n]$. In particular, 
\begin{equation*}
    S^{\sigma,\mu,W_2} =S^{\sigma,\smop(\sigma,\mu,P),W_2} \quad ~\textrm{iff}~\quad 
W_2^2(\sigma,\mu)=\sum_{i=1}^nW_2^2(\sigma^{\theta_i},\mu^{\theta_i}).
\end{equation*}    
Moreover, the registrations eliminate the effects of translation in the following sense
\begin{equation}\label{eq:registr-equal-mean}
    E({S^{\sigma,\mu,W_2}}_{\sharp}\sigma) = E({S^{\sigma,\smop(\sigma,\mu,P), W_2}}_{\sharp}\sigma) = E(\mu).
\end{equation}
\end{proposition}

\begin{proof}
Similar to the proof of \Cref{prop:equal_mean}, the minimization problems \eqref{eq:S_registr} are quadratic in terms of parameters $a$ and $b$, where $S(x)=ax+b$. By taking the partial derivatives and setting them to zero, checking the Hessian matrix, together with the assumption of the proposition, we obtain existence and uniqueness of the minimizers. The equalities then follow from direct computations. See \Cref{pf:key-affine_registration} for details.
\end{proof}
\begin{remark}
With \eqref{eq:WassDist_relation}, \Cref{prop:affine_approx} implies 
\begin{equation*}
   S^{\sigma,\mu,W_2} =S^{\sigma,\smop(\sigma,\mu,P),W_2} \quad ~\textrm{iff}~\quad 
W_2^2(\sigma,\mu)=W_2^2(\sigma,\smop(\sigma,\mu,P)).  
\end{equation*}
\Cref{prop:one_step_recovery} and \Cref{prop:one_step_recovery_compatible} show special cases where maps $S^{\sigma,\mu,W_2} $ and $S^{\sigma,\smop(\sigma,\mu,P),W_2}$ are equal. However, for the case of registration with only shifts, $S_t^{\sigma,\mu}=S_t^{\sigma,\smop(\sigma,\mu,P)}$ always holds, see \eqref{equal-shift-registration}.
\end{remark}
We further establish a connection between the Wasserstein distance and sliced-Wasserstein distance by comparing the registration maps in \eqref{eq:S_registr}. This insight holds the potential to enhance our understanding of the distinction between sliced-Wasserstein and Wasserstein flows, as demonstrated in a special case by Bonet et al. \cite[p.7, Eq. 19]{bonet2022efficient}.

\begin{corollary}
  Let $\bP$ be a random variable corresponding to the Haar probability measure $u_n$ on the orthogonal group $O(n)$. For fixed $\sigma\in \W_{2,ac}(\R^n)$ and $\mu\in \W_2(\R^n)$, we have
    \begin{align*} 
\|S^{\sigma,\mu,W_2}-\mathbb{E}S^{\sigma,\smop(\sigma,\mu,\bP), W_2}\|_{\sigma} & =  \mathbb{E}\|S^{\sigma,\mu,W_2}-S^{\sigma,\smop(\sigma,\mu,\bP), W_2}\|_{\sigma}  \\
&=\frac{W_2^2(\sigma,\mu)-nSW^2_2(\sigma,\mu)}{2\sqrt{M_2(\sigma)-\|E(\sigma)\|^2}}\geq 0.
    \end{align*}  
\end{corollary}

\begin{proof}
   By \eqref{exp-diffa}, \eqref{exp-diffb}, \eqref{eq:S-Sp} and calculations in \Cref{keyproof:prop:affine_approx}, both equations reduce to the following:
   \begin{equation} \int_{O(n)}\sum_{i=1}^nW_2^2(\sigma^{\theta_i},\mu^{\theta_i})du_n(P) = nSW_2^2(\sigma,\mu), \label{nslice-vs-SW}
   \end{equation}
which can be observed by an  explicit geometric construction of the Haar measure on $O(n)$ (see e.g., \cite[p.19]{meckes_2019}).    \end{proof}

In light of the inherent connection between $\sum_{i=1}^nW_2^2(\sigma^{\theta_i},\mu^{\theta_i})$ and $SW_2^2(\sigma,\mu)$ as shown in \eqref{nslice-vs-SW}, and the  minimization problem \eqref{eq:semi-min-problem} associated with $T_{\sigma,\mu;P}$, we demonstrate a similar connection between the registration maps associated with two distinct registration problems concerning the $W_2$ and $SW_2$ distances respectively:
\begin{corollary}\label{Cor:W2vsSW2}
Given the same assumptions as in \Cref{prop:affine_approx}, we have
    \begin{equation*}
    \|S^{\sigma,\smop(\sigma,\mu,P),W_2}-S^{\sigma,\mu,SW_2}\|_{\sigma}= \frac{nSW_2^2(\sigma,\mu)-\sum_{i=1}^n W_2^2(\sigma^{\theta^i},\mu^{\theta_i})}{\sqrt{M_2(\sigma)-\|E(\sigma)\|^2}}.
    \end{equation*}
Moreover, let $\bP$ be a random variable corresponding to the Haar probability measure $u_n$ on the orthogonal group $O(n)$, then
    \begin{equation*}
        \mathbb{E}S^{\sigma,\smop(\sigma,\mu,\bP),W_2} = S^{\sigma,\mu,SW_2}.
    \end{equation*}
\end{corollary}
\begin{proof}
    The proof follows directly from \Cref{keyproof:prop:affine_approx} and \Cref{keycor:prop:affine_approx} and similar calculations as in \eqref{exp-diffa} and \eqref{exp-diffb}.
\end{proof}
In simpler terms, the map that optimally aligns $\sigma$ with $\smop(\sigma,\mu,\bP)$ considering shifts and scalings in the $W_2$ distance is, on average, the same map that optimally aligns $\sigma$ with $\mu$ in the $SW_2$ distance.

\begin{remark}
   Similar to the registration problems \eqref{eq:S_registr} and \eqref{eq:St_registr}, one can study registration in terms of maps $\affset(P):= \{x \mapsto P\Lambda P^t x+b: \Lambda \textrm{~positive, diagonal matrix~}, b\in \R^n\}$. We provide a summary in \Cref{append:registr-compatible}.
\end{remark}

\section{Conclusion} \label{sec:conclusion}

Optimal transport-based slice-matching schemes benefit from closed-form formulations, computational efficiency and convergence guarantees. In the present paper, we are interested in the approximation power of one step of such schemes, both on the level of measures and maps. This can be considered as a step towards understanding to what extend slice-matching maps can serve as effective alternatives to optimal transport maps.

We investigate the exact recovery of basic transformations, such as translations and scalings, as well as the approximate recovery of perturbations of such transformations. These results are derived by studying invariance properties of an associated slice-matching operator. In addition, we explore equivariance and Lipschitz properties of the same operator, to understand how it incorporates actions of basic transformations on the target measure, as well as perturbations on the slicing directions. 
 
We provide a quantitative perspective on how slice-matching procedures encode special affine transformations in their approximations through the study of basic registration problems. These registration problems potentially also offer insights into the relationship between Wasserstein and sliced-Wasserstein flows, which is an interesting problem for future research. 

\section*{Acknowledgements}
CM is supported by NSF award DMS-2306064 and by a seed grant from the School of Data Science and Society at UNC. SL appreciates helpful discussions with Dr. Hengrong Du regarding \Cref{example_singleSlice} and \Cref{proofTranslationConv}.

\section*{References}
\printbibliography[title={},heading=none]

\appendix

\section{Proofs for \Cref{sec:slice-operator}}\label{appendix:proof_rigid_transformations}

\subsection{Key facts for proof of  \Cref{recovery_iff}}

 \begin{proposition}\label{prop:intersect_compatible}
 Let $\mathcal{D}(\R^n)$ be the set of differentiable vector fields from $\R^n$ to $\R^n$.

\begin{equation*}%\label{eq:intersect_compatible}
    \Big(\bigcap_{P \in O(n)}\mathfrak{S}(P)\Big) \cap \mathcal{D}(\R^n)= \{x \mapsto a\bx +\bb: a>0 \text{ and } \bb \in \R^n \}.
\end{equation*}    
\end{proposition}

\begin{proof}
    For the proof, we need to show that a differentiable vector field $S \in \bigcap_{P \in O(n)}\mathfrak{S}(P)$ is an isotropic scaling with translation. Choose $P \in O(n)$ and write $S(x) = \sum_{i=1}^n f_i^P(x\cdot \theta_{i})\theta_{i}$ with $P=[\theta_1,\ldots,\theta_n]$. Note that using the standard basis, we can also write $S(x)=\sum_{i=1}^n g_i(x_i)e_{i}$.
Computing the Jacobian of $S$ with respect to the two basis representations, we obtain
\begin{equation*}
    \begin{bmatrix}
        g_1'(x_1)&\\&\ddots \\& \quad &g_n'(x_n)
    \end{bmatrix} = P  \begin{bmatrix}
    {f_1^P}^'(\bx\cdot \bth_1)& &&\\
    & {f_2^P}^'(\bx\cdot \bth_2)&&\\
    && \ddots &\\
    &&& \ {f_n^P}'(\bx\cdot \bth_n)
    \end{bmatrix}P^t.
\end{equation*} 
Hence the two diagonal matrices above have the same diagonal entries, allowing for a possible reordering of the entries. Without loss of generality, we assume that $g_i^{\prime}(x_i)={f_i^P}^'(\bx\cdot \bth_i), i= 1,...,n$ by possibly performing a column permutation of $P$ and renaming $f_i^P$'s. Choosing an orthogonal matrix $P$ such that one of its column $\theta_i$ with all entries being non-zero, one can immediately derive that the diagonal entries $g_i^{\prime}(x_i)$'s are the same for any fixed $x$. In summary,
\begin{equation*}
    g_1^{\prime}(x_1)= \cdots = g_n^{\prime}(x_n) = {f_1^P}^{\prime}(\bx\cdot \bth_1) = \cdots = {f_n^P}^{\prime}(\bx\cdot \bth_n) = a_\bx,
\end{equation*}
where $a_\bx$ is a constant depending on $\bx = [x_1, \cdots,x_n]^t\in\R^n$. Since the diagonal element $g_i^{\prime}(x_i)$ only depends on $x_i$, it follows that $a_{\bx}$ is a constant independent of $\bx$. Hence $S(\bx)=a\bx +\bb$ for some $a>0, \bb\in\R^n$.
\end{proof}

\begin{remark}
    In general, if $T\in \bigcap_{P \in O(n)}\mathfrak{S}(P)$ is differentiable on an open set $\Omega\subseteq{\R^n}$, the $T|_{\Omega}: \Omega\rightarrow \R^n$ is an isotropic scaling with translation. In particular, $\bigcap_{P \in O(n)}\mathfrak{S}(P)$ include some piecewise isotropic scalings with translations.
\end{remark}

\subsection{Proof of \Cref{prop:Lipbounds}}\label{append:proof_Lipbounds}
We need the following proposition to derive the proof of \Cref{prop:Lipbounds}:

\begin{proposition}\label{prop:anglebound}
Consider two angles $\theta,\nu \in S^{n-1}$, and assume that $T_{\sigma^{\nu}}^{\mu^{\nu}}$ is $L$-Lipschitz for all $\nu$, i.e.\ there exists $L>0$ such that $|T_{\sigma^{\nu}}^{\mu^{\nu}}(x)-T_{\sigma^{\nu}}^{\mu^{\nu}}(y)| \leq L |x-y|$ for $x,y\in \R$ and $\nu\in S^{n-1}$, then
\begin{equation*}
    \|T_{{\sigma}^{\theta}}^{\mu^{\theta}}\circ \mathcal{P}_{\theta} - T_{{\sigma}^{\nu}}^{\mu^{\nu}}\circ \mathcal{P}_{\nu}\|_{\sigma} \leq (2L+1)C\|\theta - \nu\|_2,
\end{equation*}
where $C$ is the max over the second moments of $\sigma$ resp.\ $\mu$.
\end{proposition}
\begin{proof}
\begin{align*}
        \|T_{{\sigma}^{\theta}}^{\mu^{\theta}}\circ \mathcal{P}_{\theta} - T_{{\sigma}^{\nu}}^{\mu^{\nu}}\circ \mathcal{P}_{\nu}\|_{\sigma} \leq 
        \| T_{{\sigma}^{\theta}}^{\mu^{\theta}}\circ \mathcal{P}_{\theta} - T_{{\sigma}^{\theta}}^{\mu^{\nu}}\circ \mathcal{P}_{\theta}\|_{\sigma} + \|T_{{\sigma}^{\theta}}^{\mu^{\nu}}\circ \mathcal{P}_{\theta} - T_{{\sigma}^{\nu}}^{\mu^{\nu}}\circ \mathcal{P}_{\nu}\|_{\sigma} = (\diamond).
\end{align*}
We bound these separately.
\begin{align*}
    \| T_{{\sigma}^{\theta}}^{\mu^{\theta}}\circ \mathcal{P}_{\theta} - T_{{\sigma}^{\theta}}^{\mu^{\nu}}\circ \mathcal{P}_{\theta}\|_{\sigma} & = 
    \|T_{{\sigma}^{\theta}}^{\mu^{\theta}} - T_{{\sigma}^{\theta}}^{\mu^{\nu}}\|_{\sigma^{\theta}}
    %= \operatorname{LOT}_{\sigma^{\theta}}\left(\mu^{\theta}, \mu^{\nu}\right)
    = W_2\left(\mu^{\theta}, \mu^{\nu}\right)\\
    &\leq \|\mathcal{P}_{\theta}-\mathcal{P}_{\nu}\|_{\mu} = \left(\int_{\mathbb{R}^n}|\mathcal{P}_{\theta}(x)-\mathcal{P}_{\nu}(x)|^2\, d\mu(x)\right)^{1/2}\\
    & = \left(\int_{\mathbb{R}^n}|(\theta - \nu)\cdot x|^2\, d\mu(x)\right)^{1/2}\\
    &\leq \|\theta - \nu\|_2 \left(\int_{\mathbb{R}^n}\|x\|^2\, d\mu(x)\right)^{1/2}\\
    & \leq C\|\theta - \nu\|_2,
\end{align*}
with $C$ max of the second moments, which is bounded by assumption. Now for the second part, note that on $\mathbb{R}$ we have $T_{\sigma^{\theta}}^{\mu^{\nu}} = T_{\sigma^{\nu}}^{\mu^{\nu}} \circ T_{\sigma^{\theta}}^{\sigma^{\nu}}$
\begin{align*}
   \|T_{{\sigma}^{\theta}}^{\mu^{\nu}}\circ \mathcal{P}_{\theta} - T_{{\sigma}^{\nu}}^{\mu^{\nu}}\circ \mathcal{P}_{\nu}\|_{\sigma}  =  
   \left(\int_{\mathbb{R}^n} |T_{\sigma^{\nu}}^{\mu^{\nu}}(T_{\sigma^{\theta}}^{\sigma^{\nu}}( \mathcal{P}_{\theta}(x))) - T_{{\sigma}^{\nu}}^{\mu^{\nu}}(\mathcal{P}_{\nu}(x))|^2\, d\sigma(x)\right)^{1/2} = (\star)
\end{align*}
Since $T_{\sigma^{\nu}}^{\mu^{\nu}}$ is $L$-Lipschitz, we get

\begin{align*}
    (\star) & \leq L\,\left(\int_{\mathbb{R}^n} | T^{\sigma^{\nu}}_{\sigma^{\theta}}(\mathcal{P}_{\theta}(x)) -  \mathcal{P}_{\nu}(x)|^2\, d\sigma(x)\right)^{1/2} 
    = L\|T^{\sigma^{\nu}}_{\sigma^{\theta}}\circ \mathcal{P}_{\theta} -\mathcal{P}_{\nu} \|_{\sigma} \\
    & \leq L\left(\|T^{\sigma^{\nu}}_{\sigma^{\theta}}\circ \mathcal{P}_{\theta} - \mathcal{P}_{\theta}\|_{\sigma} + \|\mathcal{P}_{\theta}-\mathcal{P}_{\nu}\|_{\sigma}\right)\\
    & \leq L \left(\|T^{\sigma^{\nu}}_{\sigma^{\theta}} - \operatorname{id}\|_{\sigma^{\theta}} + C\|\theta-\nu\|_2 \right) \\
    & = L\left(W_2(\sigma^{\theta},\sigma^{\nu})+C \|\theta - \nu\|_2 \right) \\
    & \leq L\left(\|\mathcal{P}_{\theta}-\mathcal{P}_{\nu}\|_{\sigma}+C \|\theta - \nu\|_2 \right) \\
    & \leq 2LC\|\theta - \nu\|_2
\end{align*}
This implies
\begin{equation*}
    (\diamond) \leq (2L+1)C\|\theta - \nu\|_2. \qedhere
\end{equation*}
\end{proof}

%%%%%%% PROOF OF STABILITY THEOREM FOR 1 STEP
\begin{proof}[Proof of \Cref{prop:Lipbounds}]
Based on \eqref{eq:compatibleSets}, we let $T_{\sigma,\mu;P} = PD\circ P^t$ where $D(x)=[T_{\sigma^{\theta_1}}^{\mu^{\theta_1}}(x_1), 
     T_{\sigma^{\theta_2}}^{\mu^{\theta_2}}(x_2),
   \cdots,
    T_{\sigma^{\theta_n}}^{\mu^{\theta_n}}(x_n)]^t$ for $x \in \R^n$ and 
 $P = [\theta_1,\ldots,\theta_n]$. Similarly, we let and $T_{\sigma,\mu;Q} = Q\widetilde{D}\circ Q^t$,  with $Q = [\nu_1,\ldots,\nu_n]$. We continue with deriving the bound:
%%%%%%%%%%%% BOUNDS %%%%%%%%%%%%%%%%%%%%%%%%%
\begin{align*}
    \|T_{\sigma,\mu;P}-T_{\sigma,\mu;Q}\|_{\sigma} &= \|PDP^t-Q\widetilde{D}Q^t\|_{\sigma}\\
    &\leq \|PDP^t - P\widetilde{D}Q^t\|_{\sigma} + \|P\widetilde{D}Q^t-Q\widetilde{D}Q^t\|_{\sigma} \\
    & = (1) + (2).
\end{align*}
We bound the two terms seperately. For (1), using \Cref{prop:anglebound}, we get
\begin{align*}
    \|PDP^t - P\widetilde{D}Q^t\|_{\sigma}^2 
    & = \int_{\R^n} \|D(P^tx) - \widetilde{D}(Q^tx)\|_2^2\,d\sigma(x) \\
    & = \sum_{i=1}^n \int_{\R^n} |T_{\sigma^{\theta_i}}^{\mu^{\theta_i}}((P^tx)_i)-T_{\sigma^{\nu_i}}^{\mu^{\nu_i}}((Q^tx)_i)|^2\, d\sigma(x) \\
    & = \sum_{i=1}^n \|T_{\sigma^{\theta_i}}^{\mu^{\theta_i}}\circ \mathcal{P}_{\theta_i} - T_{\sigma^{\nu_i}}^{\mu^{\nu_i}}\circ \mathcal{P}_{\nu_i} \|_{\sigma}^2 \\
    & \leq ((2L+1)C)^2 \sum_{i=1}^n \|\theta_i - \nu_i\|_2^2 \\
    & = ((2L+1)C)^2 \|P - Q\|_F^2.
\end{align*}
For (2) we get
\begin{align*}
    \|P\widetilde{D}Q^t-Q\widetilde{D}Q^t\|_{\sigma}^2 & = \int_{\R^n} \|(P-Q)\widetilde{D}(Q^tx)\|_2^2 \, d\sigma(x)  \\
   & \leq \|P-Q\|_2^2\int_{\R^n} \|\widetilde{D}(Q^tx)\|_2^2 \, d\sigma(x) \\
    & \leq \|P-Q\|_2^2 \, L^2 \int_{\R^n} \|Q^tx\|_2^2\, d\sigma(x) 
     \leq \|P-Q\|_2^2 \, L^2 C^2   \\
     & \leq \|P-Q\|_F^2 \, L^2 C^2
\end{align*}
Combining (1) and (2) gives the final bound.
\end{proof}

\section{Proofs for \Cref{sec:registration}}

\begin{proof}[Proof of  \Cref{prop:affine_approx}]\label{pf:key-affine_registration}
Let $ S^{\sigma,\mu,W_2}(x) =  a^{W_2} x+  b^{W_2}$ and $S^{\sigma,\smop(\sigma,\mu,P), W_2}(x)=  \widetilde{a}^{W_2}x+ \widetilde{b}^{W_2}$ be the critical functions for the associated minimization problem \eqref{eq:S_registr}. By \Cref{keyproof:prop:affine_approx} and \Cref{keycor:prop:affine_approx}, we have
 \begin{align}
\widetilde{a}^{W_2} - a^{W_2} &=\frac{W_2^2(\sigma,\mu)-\sum_{i=1}^nW^2_2(\sigma^{\theta_i},\mu^{\theta_i})}{2(M_2(\sigma)-\|E(\sigma)\|^2)},\label{exp-diffa}\\
\widetilde{b}^{W_2}-b^{W_2} &= -(\widetilde{a}^{W_2} - a^{W_2})E(\sigma),\label{exp-diffb}
   \end{align}
and the norm bound $ \|S^{\sigma,\mu,W_2}-S^{\sigma,\smop(\sigma,\mu,P), W_2}\|_{\sigma}$ in \eqref{eq:S-Sp}
can be obtained via direct computation and the fact that the RHS is non-negative, see \Cref{lm:WassSlice_compare}. It is left to show that these critical functions are indeed the minimizers by verifying 
\begin{enumerate}
    \item $\widetilde{a}^{W_2}\geq a^{W_2} >0$, see \Cref{lm:WassSlice_compare}, \Cref{lm:pos_denom}, and \Cref{lm:nonneg_numer}.
    \item The Hessian associated $H(a,b)$ with both the minimization problems are positive definite by a direct calculation and \Cref{lm:pos_denom}, where
    \begin{equation*}
    H(a,b)=   2 \begin{bmatrix}
           M_2(\sigma) & (E(\sigma))^t\\
           E(\sigma)    & I_{n-1}
        \end{bmatrix}.
    \end{equation*}
Here $I_{n-1}$ denotes the identity matrix of size $(n-1)\times (n-1)$ .
\end{enumerate}    
The equality concerning the means follows from \Cref{cor:registr-mean}.    
\end{proof}
\begin{proposition}\label{keyproof:prop:affine_approx}
    Let $S^{\sigma,\eta,W_2}$ and $S^{\sigma,\eta, SW_2}$ correspond to the critical points of the minimization problems in \eqref{eq:S_registr} and \eqref{eq:SP_registr}, respectively. Then the corresponding parameters 
    satisfy
    \begin{align*}
    a^{W_2} &=  \frac{\frac{1}{2}(M_2(\eta)+M_2(\sigma)-W_2^2(\sigma,\eta))-E(\sigma)\cdot E(\eta)}{M_2(\sigma)-\|E(\sigma)\|^2},\\
    b^{W_2} & = E(\eta)-a^{W_2}E(\sigma),\\
     a^{SW_2} &=  \frac{\frac{1}{2}(M_2(\eta)+M_2(\sigma)-nSW_2^2(\sigma,\eta))-E(\sigma)\cdot E(\eta)}{M_2(\sigma)-\|E(\sigma)\|^2},\\
    b^{SW_2} & = E(\eta)-a^{SW_2}E(\sigma),  
\end{align*}
where $S^{\sigma,\eta,W_2}(x) =  a^{W_2} x+  b^{W_2}$ and $S^{\sigma,\eta,SW_2}(x) =  a^{SW_2} x+  b^{SW_2}$.
\end{proposition}

\begin{proof}
Given $\sigma\in \W_{2,ac}(\R^n)$ and $\eta\in \W_2(\R^n)$, let $M_2(\sigma)= \int \|x\|^2d\sigma(x)$ (similarly define $M_2(\eta)$), $E(\sigma)= \int xd\sigma(x)$ (similarly define $E(\eta)$). For $S(x)=ax+b$, by the changes of variables formula and the fact that $T_{\sigma}^\eta = T_{S_{\sharp}\sigma}^{\eta}\circ S$, we have
\begin{align*}
 W_2^2(S_{\sharp}\sigma,\eta)& =  \|T_{\sigma}^{\eta}-(ax+b)\|^2_{\sigma} = M_2(\eta)+a^2M_2(\sigma)+2ab\cdot E(\sigma)\\&-2a\int T_{\sigma}^{\eta}(x)\cdot x d\sigma(x)-E(\sigma)\cdot E(\eta)+\|b\|^2-2E(\eta)\cdot b.
\end{align*}
Taking the partial derivatives gives
\begin{align*}
    \frac{\partial}{\partial a} &= 2a M_2(\sigma)+2b\cdot E(\sigma)-2\int T_{\sigma}^{\eta}(x)\cdot x d\sigma(x), \\
    \frac{\partial}{\partial b} & = 2b+2aE(\sigma)-2E(\eta).
\end{align*}
Setting the above equations to zero and with the observation that $\int T_{\sigma}^{\eta}(x)\cdot x d\sigma(x) = \frac{1}{2}(M_2(\eta)+M_2(\sigma)-W_2^2(\sigma,\eta))$, we get the the desired formulas for $a^{W_2}$ and $b^{W_2}$.
Similarly, 
\begin{align*}
 SW_2^2(S_{\sharp}\sigma,\eta) &=  \int_{S^{n-1}}W_2^2((S_{\sharp}\sigma)^{\theta}, \eta^{\theta})du(\theta) \\&= \int_{S^{n-1}}\int_{\R} |T_{\sigma^{\theta}}^{\eta^\theta}(t)-(at+b\cdot\theta)|^2dt du(\theta)\\& = \frac{1} {n}\Big(M_2(\eta)+a^2M_2(\sigma)+2ab\cdot E(\sigma)\\&-2na\int_{S^{n-1}}\int_{\R} tT_{\sigma^{\theta}}^{\eta^\theta}(t)d\sigma^{\theta}(t) du(\theta)-E(\sigma)\cdot E(\eta)+\|b\|^2-2E(\eta)\cdot b\Big).
\end{align*}
Taking the partial derivatives gives
\begin{align*}
    \frac{\partial}{\partial a} &=\frac{1}{n} \Big(2a M_2(\sigma)+2b\cdot E(\sigma)-2n\int_{S^{n-1}}\int_{\R} tT_{\sigma^{\theta}}^{\eta^\theta}(t)dt du(\theta)\Big), \\
    \frac{\partial}{\partial b} & = \frac{1}{n} \Big(2b+2aE(\sigma)-2E(\eta)\Big).
\end{align*}
Setting the above equations to zero and with the observation that $\int_{S^{n-1}}\int_{\R} tT_{\sigma^{\theta}}^{\eta^\theta}(t)d\sigma^{\theta}(t) du(\theta) = \frac{1}{2n}(M_2(\sigma)+M_2(\eta)-nSW_2^2(\sigma,\eta))$ , we get the desired formulas for $a^{SW_2}$ and $b^{SW_2}$. We provide computational details in \Cref{append:details}.
\end{proof}

\begin{corollary}\label{cor:registr-mean}
Given the same assumptions as in \Cref{keyproof:prop:affine_approx}, for $D=W_2~\textrm{or}~ SW_2$ 
    \begin{equation}
        E({S^{\sigma,\eta,D}}_\sharp\sigma) = E(\eta).
    \end{equation}
\end{corollary}
\begin{proof}
    Upon direct calculation,  we have $E({S^{\sigma,\eta,D}}_\sharp\sigma = a^DE(\sigma)+b^D$,  where $a^D, b^D$ are as in  \Cref{keyproof:prop:affine_approx}. The conclusion can be derived from the expressions for $b^{W_2}$ and $b^{SW_2}$.
\end{proof}

\begin{corollary}\label{keycor:prop:affine_approx}
Let $\eta = \smop(\sigma,\mu,P)$ in \Cref{keyproof:prop:affine_approx}. Then the parameters corresponding to   $S^{\sigma,\smop(\sigma,\mu,P), W_2}$ and
 $S^{\sigma,\smop(\sigma,\mu,P), SW_2}$ satisfy 
\begin{align*}
    \widetilde{a}^{W_2} & =  \frac{\frac{1}{2}(M_2(\mu)+M_2(\sigma)-\sum_{i=1}^nW_2^2(\sigma^{\theta_i},\mu^{\theta_i}))-E(\sigma)\cdot E(\mu)}{M_2(\sigma)-\|E(\sigma)\|^2},\\
    \widetilde{b}^{W_2} & = E(\mu)-\widetilde{a}^{W_2}E(\sigma),\\
        \widetilde{a}^{SW_2} &=  \frac{\frac{1}{2}(M_2(\mu)+M_2(\sigma)-nSW_2^2(\sigma,\smop(\sigma,\mu,P)))-E(\sigma)\cdot E(\mu)}{M_2(\sigma)-\|E(\sigma)\|^2},\\
   \widetilde{b}^{SW_2} & = E(\mu)-a^{SW_2}E(\sigma),
\end{align*}
where $S^{\sigma,\smop(\sigma,\mu,P), W_2}(x)=  \widetilde{a}^{W_2}x+ \widetilde{b}^{W_2}$  and $S^{\sigma,\smop(\sigma,\mu,P), SW_2}(x)=  \widetilde{a}^{SW_2}x+ \widetilde{b}^{SW_2}$. 
\end{corollary}
\begin{proof}
The above formulas follows directly from \Cref{keyproof:prop:affine_approx}, the fact that $\smop(\sigma,\mu,P)$ and $\mu$ have the same mean
(see \eqref{eq:equalmean}), and the formula \eqref{eq:WassDist_relation} for $W_2^2(\sigma, \smop(\sigma,\mu,P))$.
\end{proof}

\begin{proposition}\label{append:registr-compatible}
Let $$\affset(P):= \{x\mapsto P\Lambda P^t x+b: \Lambda \textrm{~positive, diagonal matrix~}, b\in \R^n\}.$$ Consider the minimization problem
\begin{align}
    S_P^{\sigma,\eta} &:= \argmin_{S_P\in \affset(P)} \|{S_P}-T_\sigma^\eta\|_{\sigma}.\label{eq:SP_registr}
\end{align}
Let $S^{\sigma,\mu}_{P}$ and $S^{\sigma,\smop(\sigma,\mu;P)}_{P}$ be the minimizers of \eqref{eq:SP_registr} with $\eta=\mu$ and $\eta = \smop(\sigma,\mu,P)$, respectively. We denote the diagonal entries of the corresponding $\Lambda$ by $a_i$ and $\widetilde{a}_i$, respectively. Similar notation holds for $b_i$ and $\widetilde{b}_i$. Then 
    \begin{align*}
        \widetilde{a}_i - a_i &= \frac{\int|\theta_i\cdot (T_{\sigma}^{\mu}(x)-x)|^2d\sigma(x)- W_2^2(\sigma^{\theta_i},\mu^{\theta_i})}{2(M_2^{\sigma^{\theta_i}}- (E^{\theta_i})^2)}\geq 0,\\
         \widetilde{b}_i - b_i &= - \sum_{i=1}^n\theta_i E^{\sigma^{\theta_i}} (  \widetilde{a}_i - a_i).
    \end{align*}
\end{proposition} 
\begin{proof}
    The proof uses similar arguments in \Cref{keyproof:prop:affine_approx} and \Cref{keycor:prop:affine_approx} except the partial derivatives are with respect to $a_i$ and $\widetilde{a_i}$ instead of $a$ and $\widetilde{a}$. Note that following these arguments, we use the equations presented in \Cref{lm:SWeq}.
\end{proof}

\section{Other technical details}\label{append:details}
\begin{lemma}\label{lm:equal_mean}
Let $\sigma\in \W_{2,ac}(\R^n)$ and $\eta, \mu\in \W_2(\R^n)$. Then we get
    \begin{align*} %\label{eq:equal_mean}     
        &E({(T_{\sigma,\mu;P})_{\sharp}\sigma})=\int T_{\sigma,\mu;P}(x) d\sigma(x) = \int yd\mu(y)=E(\mu),\\
        &\int T_{\sigma}^{\eta}(x)\cdot x d\sigma(x) = \frac{1}{2}(M_2(\eta)+M_2(\sigma)-W_2^2(\sigma,\eta))\\
       & M_2({(T_{\sigma,\mu;P})_{\sharp}\sigma})=\int \| T_{\sigma,\mu;P}(x)\|^2 d\sigma(x) = M_2(\mu)%\label{equal2ndMoment}
    \end{align*}
\end{lemma}
\begin{proof}
By the change of variables formula, we have
   \begin{align*}
        \int T_{\sigma,\mu;P}(x) d\sigma(x)& = \sum_{i=1}^n \theta_i\int T_{\sigma^{\theta_i}}^{\mu^{\theta_i}}(x\cdot \theta_i)d\sigma(x)
        = \sum_{i=1}^n \theta_i\int  T_{\sigma^{\theta_i}}^{\mu^{\theta_i}}(t)d\sigma^{\theta_i}(t)\\
        &= \sum_{i=1}^n \theta_i\int zd\mu^{\theta_i}(z)
        = \sum_{i=1}^n \theta_i\int y\cdot \theta_i d\mu(y)\\
        & = \int y d\mu(y),
   \end{align*} 
\begin{align*}
    \int T_{\sigma}^{\eta}(x)\cdot x d\sigma(x) &= \frac{1}{2}\Big(\int \|T_{\sigma}^{\eta}(x)\|^2d\sigma(x)\\&+\int \|x\|^2d\sigma(x)- \int \|T_{\sigma}^{\eta}(x)-x\|^2d\sigma(x)\Big)\\&=\frac{1}{2}(M_2(\eta)+M_2(\sigma)-W_2^2(\sigma,\eta)),
\end{align*}
\begin{align*}
    \int \| T_{\sigma,\mu;P}(x)\|^2 d\sigma(x) & =\int \sum_{i=1}^n |T_{\sigma^{\theta_i}}^{\mu^{\theta_i}}(x\cdot \theta)|^2d\sigma(x)
     = \sum_{i=1}^n \int|T_{\sigma^{\theta_i}}^{\mu^{\theta_i}}(t)|^2d\sigma^{\theta_i}(t)\\
    & =\sum_{i=1}^n \int |w|^2d\mu^{\theta_i}(w) = \sum_{i=1}^n\int|y\cdot \theta_i|^2d\mu(y)\\
    & = \int \|y\|^2d\mu(y)= M_2(\mu),
\end{align*}
where the last steps make use of the fact that $P=[\theta_1,\cdots,\theta_n]$ is an orthogonal matrix.
\end{proof}

\begin{lemma}\label{lm:SWeq}
Let $\sigma\in \W_{2,ac}(\R^n)$, $\eta\in \W_2(\R^n)$, and $b\in \R^n$. Then
    \begin{align}
       &\int_{S^{n-1}}\int_{\R} tT_{\sigma^{\theta}}^{\eta^\theta}(t)d\sigma^{\theta}(t) du(\theta) = \frac{M_2(\sigma)+M_2(\eta)-nSW_2^2(\sigma,\eta)}{2n}\label{t-mixed-terms}\\
   &\int_{S^{n-1}}\int_{\R} t^2d\sigma^{\theta}(t)du(\theta) = \frac{M_2(\sigma) }{n}\label{M2sigma}\\
      &\int_{S^{n-1}}\int_{\R} |T_{\sigma^{\theta}}^{\eta^\theta}(t)|^2d\sigma^{\theta}(t)du(\theta) = \frac{M_2(\eta)}{n}\label{M2eta}\\
   &\int_{S^{n-1}}\int_{\R} (b\cdot \theta) t d\sigma^{\theta}(t)du(\theta)  = \frac{E(\sigma)\cdot b}{n}\label{Ebsigma}\\
     &\int_{S^{n-1}}\int_{\R} (b\cdot\theta) T_{\sigma^{\theta}}^{\eta^\theta}(t)td\sigma^{\theta}(t)du(\theta)  = \frac{E(\eta)\cdot b}{n}\label{Ebeta}
         \end{align}
\end{lemma}
\begin{proof}
 We note that \eqref{M2sigma} and \eqref{M2eta} are analogous by the change of variables formula, so are \eqref{Ebsigma} and \eqref{Ebeta}. We will first show  \eqref{M2sigma}.
 \begin{align*}
  \int_{S^{n-1}}\int_{\R} t^2d\sigma^{\theta}(t)du(\theta) & = \int_{S^{n-1}}\int_{\R^n}|x\cdot \theta|^2 d\sigma(x)du(\theta) \\
   & \hspace{{-.4cm}}\overset{\textrm{Fubini}}{=} \int_{\R^n}\int_{S^{n-1}}|x\cdot \theta|^2 du(\theta)d\sigma(x)\\
  &= \int_{\R^n}\frac{\|x\|^2}{2}d\sigma(x)\\
  &=\frac{M_2(\sigma)}{n}.
 \end{align*}
For \eqref{Ebsigma}, we have
\begin{align*}
  \int_{S^{n-1}}\int_{\R} (b\cdot \theta) t d\sigma^{\theta}(t)du(\theta) & = \int_{S^{n-1}}b\cdot \theta \int_{\R^n}x\cdot \theta d\sigma(x)du(\theta)\\
  &= \int_{S^{n-1}}(b\cdot \theta)(E(\sigma)\cdot \theta) du(\theta)\\
&=  \int_{S^{n-1}}\frac{1}{2}\Big(|b\cdot\theta)|^2+|E(\sigma)\cdot \theta|^2-|(b-E(\sigma))\cdot\theta|^2\Big)du(\theta)\\
&= \frac{1}{2n}\Big(\|b\|^2+\|E(\sigma)\|^2- \|b-E(\sigma)\|^2\Big)\\
& = \frac{E(\sigma)\cdot b}{n}.
\end{align*}
With  \eqref{M2sigma} and \eqref{M2eta}, we have \eqref{t-mixed-terms}:
\begin{align*}
  &\int_{S^{n-1}}\int_{\R} tT_{\sigma^{\theta}}^{\eta^\theta}(t)d\sigma^{\theta}(t) du(\theta) \\& =\frac{1}{2} \int_{S^{n-1}}\int_{\R}\Big (t^2+(T_{\sigma^{\theta}}^{\eta^\theta}(t))^2- (t-T_{\sigma^{\theta}}^{\eta^\theta}(t))^2\Big)  d\sigma^{\theta}(t) du(\theta)\\
  & = \frac{1}{2n}\Big(M_2(\sigma)+M_2(\eta)- n\int_{S^{n-1}}W_2^2(\sigma^{\theta},\eta^{\theta})du(\theta)\Big)\\
  & = \frac{M_2(\sigma)+M_2(\eta)- nSW_2^2(\sigma,\eta)}{2n}.
\end{align*}

\end{proof}

\begin{lemma}\label{lm:WassSlice_compare}
Let $\sigma\in \W_{2,ac}(\R^n)$ and $\mu\in \W_2(\R^n)$ and $P = [\theta_1, \cdots, \theta_n]\in O(n)$. Then
    \begin{equation*}
         W_2^2(\sigma,\mu) \geq \sum_{i=1}^nW_2^2(\sigma^{\theta_i},\mu^{\theta_i}).
    \end{equation*}
\end{lemma}
\begin{proof}
By \cite[Proposition 5.1.3]{bonnotte13thesis}, 
\begin{equation*}
    W_2^2(\sigma^{\theta},\mu^{\theta})\leq \int |\theta\cdot x-\theta\cdot y|^2d\gamma^*(x,y),
\end{equation*}
where $\gamma^*$ is the optimal transport plan between $\sigma$ and $\mu$. Then 
\begin{align*}
 \sum_{i=1}^n W_2^2(\sigma^{\theta_i},\mu^{\theta_i})&\leq  \int \sum_{i=1}^n   |\theta_i\cdot (x-y)|^2d\gamma^*(x,y)\\
 &= \int \|x-y\|^2d\gamma^*(x,y)\\
 &= W_2^2(\sigma,\mu).\qedhere
\end{align*}
\end{proof}

\begin{lemma}\label{lm:Holder}
    Let $h:\R^n\rightarrow\R^n$ and $\sigma(\R^n)=1$. Then 
    \begin{equation*}
        \int \|h(x)\|^2d\sigma(x) \|\geq \|\int h(x)d\sigma(x)\|^2,
    \end{equation*}
    where equality holds if and only if $h(x)=v$  $\sigma$-a.e. for some $v\in \R^n$.
\end{lemma}
\begin{proof}
    Let $h(x)= [h_1(x), \cdots, h_n(x)]^t$. By H\"{o}lder's inequality, 
    \begin{align*}
  \int |h_i(x)|d\sigma(x&)\leq \left(\int |h_i(x)|^2d\sigma(x)\right)^{1/2}\left(\int 1^2 d\sigma(x)\right)^{1/2}     \\
  &= \left(\int |h_i(x)|^2d\sigma(x)\right)^{1/2}.
    \end{align*}
Squaring the above inequality and summing over $i$ gives the desired inequality. Observe that equality holds if and only if $h_i(x)=v_i$ for some constant $v_i\in \R$.
\end{proof}

\begin{lemma}\label{lm:pos_denom}
    Let $\sigma\in \W_{2,ac}(\R^n)$, and $ M_2(\sigma), E(\sigma)$ be defined as in \Cref{keyproof:prop:affine_approx}. Then 
    \begin{equation*}
        M_2(\sigma)- \|E(\sigma)\|^2>0.
    \end{equation*}
\end{lemma}
\begin{proof}
 Since $x$ is not a constant vector $\sigma$-a.e. ($\sigma\in \W_{2,ac}(\R^n)$), it follows from \Cref{lm:Holder} with $h(x)=x$ that
    \begin{equation*}
        \int \|x\|^2 d\sigma(x) > \|\int x d\sigma(x)\|^2.\qedhere
    \end{equation*}
\end{proof}

\begin{lemma}\label{lm:nonneg_numer}
Let $\sigma\in 
\W_{2,ac}(\R^n), \mu \in \W_2(\R^n)$  and $\phi$ be a convex function such that $\grad \phi = T_{\sigma}^{\mu}$ given by Brenier's theorem (see e.g., \cite[Theorem 1.48]{santambrogio2015optimal}). If $\phi$ is differentiable at $E(\sigma)$, where $E(\sigma) = \int xd\sigma(x)$, then 
\begin{equation}\label{eq:monotoneOT}
    \int T_{\sigma}^{\mu}(x)\cdot x d\sigma(x)- \Big(\int xd\sigma(x)\Big)\cdot \Big(\int T_{\sigma}^{\mu}(x)d\sigma(x)\Big)\geq 0.
\end{equation}
\end{lemma}
\begin{proof}
Let $A = \{x\in \R^n: \phi \textrm{~is differentiable at ~} x\}$.
Since $\phi$ is $\sigma$-a.e. differentiable, we have $\sigma(A)=1$. Then it follows from the convexity of $\phi$ that
\begin{equation*}%\label{eq:pos_gradphi}
    (\grad\phi(x)-\grad\phi(E(\sigma)))\cdot (x-E(\sigma))\geq 0, \quad \forall x\in A.
    \end{equation*}
Hence 
\begin{equation*}
    \int_{A} (T_{\sigma}^{\mu}(x)-T_{\sigma}^{\mu}(E(\sigma)))\cdot (x-E(\sigma)) d\sigma(x)\geq 0,
\end{equation*}
which is exactly the desired inequality \eqref{eq:monotoneOT} by a direct computation using $\sigma(A)=1$:
\begin{align*}
 &-\int T_{\sigma}^{\mu}(E(\sigma))\cdot x d\sigma(x)- \int T_{\sigma}^{\mu}(x)\cdot E(\sigma) d\sigma(x) + T_{\sigma}^{\mu}(E(\sigma))\cdot E(\sigma) \\& = -T_{\sigma}^{\mu}(E(\sigma)) \cdot E(\sigma) - \Big(\int xd\sigma(x)\Big)\cdot \Big(\int T_{\sigma}^{\mu}(x)d\sigma(x)\Big)+ T_{\sigma}^{\mu}(E(\sigma)) \cdot E(\sigma)\\
 & = \Big(\int xd\sigma(x)\Big)\cdot \Big(\int T_{\sigma}^{\mu}(x)d\sigma(x)\Big).\qedhere
\end{align*}
\end{proof}

\begin{remark}\label{alter_assum}
    The same conclusion holds if the assumption were ``$E(\sigma)$ lies in the support of $\sigma$" instead of $\phi$ being differentiable at $E(\sigma)$, which can be proved using the fact that the support of optimal transport plan is cyclically monotone.
\end{remark}

\begin{remark}\label{add_assump_strict}
    Given the assumptions in \Cref{lm:nonneg_numer}, one can show that the inequality is strict if in addition,  there exists a ball $B(x,r)$, where $x$ lies in the support of $\sigma$, such that for any $
    \lambda \in (0,1)$ and $y\in B(x,r)$
    \begin{equation*}
        \phi((1-\lambda)y+\lambda E(\sigma)) < (1-\lambda)\phi(y)+\lambda\phi(E(\sigma)),
    \end{equation*}
which guarantees that the inequality $eq:pos_gradphi$ is strict for $y$ in a set with positive measure. 
In particular, if furthermore $\phi$ in \Cref{lm:nonneg_numer} is  strictly convex, the strict inequality holds.
\end{remark}

\begin{proposition}\label{proofTranslationConv}
    Let $\sigma\in \W_{2,ac}(\R^n)$ and  $\mu = T^b_{\sharp}\sigma$ with $T^b(x)= x+b, b\neq 0\in \R^n$. Consider iteration $\sigma_{k+1}=(T_{\sigma_k,\mu; \theta_k})_\sharp \sigma_k$, with $\sigma_0=\sigma$ and where $\theta_k$ is chosen i.i.d. according to the uniform measure on $S^{n-1}$. Then 
    \begin{equation*}
        \sigma_k \xrightarrow{a.s.} \mu \quad \textrm{in~} W_2.
    \end{equation*}
\end{proposition}
  \begin{proof}
By a direct computation, $T_{\sigma_k}^{\mu}(x) =x+b_k$, where  
\begin{equation*}
        b_{k+1}= b_k - \theta_k(\theta_k\cdot b_k).
    \end{equation*}
To show $\sigma_k\rightarrow \mu$ almost surely, it suffices to show that $b_k \rightarrow 0$ almost surely.
By symmetry of $S^{n-1}$, we assume without of generality that $b_0 = [1,0,\cdots,0]^t$.
Note that $\|b_1\|^2= 1-|\theta_0\cdot b_0|^2$. 
Consider the spherical coordinates for $S^{n-1}$ with $\phi_1,\cdots \phi_{n-2}\in [0,\pi]$ and $\phi_{n-1}\in [0,2\pi]$:
\[
\begin{aligned}
&x_{1} = \cos(\varphi_{1}),\quad 
x_{2} = \sin(\varphi_{1})\cos(\varphi_{2}), \quad
x_{3} = \sin(\varphi_{1})\sin(\varphi_{2})\cos(\varphi_{3}) \\
&\hspace{6cm} \cdots  \\
&x_{n-1} = \sin(\varphi_{1})\cdots\sin(\varphi_{n-2})\cos(\varphi_{n-1}),\  x_{n} = \sin(\varphi_{1})\cdots\sin(\varphi_{n-2})\sin(\varphi_{n-1}). \\
\end{aligned}
\]
The corresponding Jacobian is $\sin^{n-2}(\varphi_1)\sin^{n-3}(\varphi_2)\cdots \sin\varphi_{n-2}$. A direct computation gives
\begin{align*}
  \mathbb{E}[|\theta_0\cdot b_0|^2] &= \frac{\int_0^{\pi} \sin^{n-2}(\varphi_1)\cos^2(\varphi_1)d\varphi_1}{\int_0^{\pi}\sin^{n-2}(\varphi_1)d\varphi_1}\\
  & = 1-\frac{\int_0^{\pi} \sin^{n}(\varphi_1)d\varphi_1}{\int_0^{\pi}\sin^{n-2}(\varphi_1)d\varphi_1}
  \\&= \rho <1.
\end{align*}
Hence $ \mathbb{E}[\|b_1\|^2] = 1-\rho \in (0,1)$. 
By symmetry and induction, one can show that 
\begin{equation*}
     \mathbb{E}[\|b_k\|^2] = (1-\rho)^k 
     \xrightarrow{k\rightarrow\infty} 0.
\end{equation*}
Since $\|b_{k+1}\|\leq \|b_k\|$, by the monotone convergence theorem, we have 
\begin{equation*}
    \mathbb{E}[\|b_k\|^2]\longrightarrow \mathbb{E}[\alpha_{\infty}^2],
\end{equation*}
where $\alpha_{\infty} = \lim \alpha_k$ and $\alpha_k = \|b_k\|$, which implies 
$\alpha_{\infty}= 0$ almost surely and hence $b_k \to 0$  almost surely.
\end{proof}

\begin{lemma}\label{lm:abs_cont}
    Let $\sigma,\mu\in \W_{2,ac}(\R^n)$. Then $(T_{\sigma,\mu;P})_\sharp \sigma\in \W_{2,ac}(\R^n)$, for any $P\in O(n)$.
\end{lemma}

\begin{proof}
    Let $P = [\theta_1,\cdots, \theta_n]$. A direct computation shows 
    \begin{align*}
      \grad T_{\sigma,\mu;P}(x) = P  \begin{bmatrix}
   (T_{\sigma^{\theta_2}}^{\mu^{\theta_i}})^{\prime}(\bx\cdot \bth_1)& &&\\
    & (T_{\sigma^{\theta_2}}^{\mu^{\theta_i}})^{\prime}(\bx\cdot \bth_2)&&\\
    && \ddots &\\
    &&& (T_{\sigma^{\theta_n}}^{\mu^{\theta_i}})^{\prime}(\bx\cdot \bth_n)
    \end{bmatrix}P^t.
    \end{align*}
Following similar arguments as in \cite[Proof of Lemma 1, p. 949]{zemel2019frechet}, it suffices to show that there exists a set $\Sigma$ such that (i) $\sigma(\R^n\setminus \Sigma)=0$ (ii) $T_{\sigma,\mu;P}|_{\Sigma}$ is injective and $ \grad T_{\sigma,\mu;P}$ is positive definite on $\Sigma$. To this end, it suffices to observe that $T_{\sigma^{\theta_i}}^{\mu^{\theta_i}}$ is injective and $(T_{\sigma^{\theta_i}}^{\mu^{\theta_i}})^{\prime}>0$  outside a set $U_i$ that is $\sigma^{\theta_i}$-negligible, i.e., $\sigma^{\theta_i}(U_i)=0$. Here we have used the fact that $T_{\sigma^{\theta_i}}^{\mu^{\theta_i}}$ exists and is unique given that $\sigma\in \mathcal{P}_{ac}(\R^n)$ (and hence $\sigma^{\theta_i}$ is absolutely continuous, see e.g., Box 2.4. in \cite[p. 82]{santambrogio2015optimal}). The fact that $M_2((T_{\sigma,\mu;P})_\sharp \sigma)$ is finite follows from \eqref{eq:equal2ndMoment} and that $M_2(\mu)<\infty$.
\end{proof}

\end{document}